\documentclass[twocolumn]{autart}    % Enable this line and disable the
\usepackage{graphicx}          % Include this line if your
                               \usepackage{amsfonts}
                               \usepackage{mathrsfs}
                               \usepackage{colortbl}
                               \usepackage{amssymb,amsmath}
                                \usepackage{epstopdf}
                               \usepackage{mathtools}
							   \usepackage{natbib}        % required for bibliography
                               \usepackage{tikz}
                               \usetikzlibrary{positioning}
                               \allowdisplaybreaks
                               \usepackage{wrapfig}

\newtheorem{Assumption}{\bf Assumption}[section]
\newtheorem{Remark}{\bf Remark}[section]

\newtheorem{Problem}{\bf Problem}
\newenvironment{Proof}{\noindent{\em Proof:\/}}{\hfill $\Box$\par}
\newtheorem{Theorem}{\bf Theorem}[section]
\newtheorem{Lemma}{\bf Lemma}[section]

\newtheorem{Definition}{\bf Definition}[section]
\newtheorem{Proposition}{\bf Proposition}[section]

\begin{document}

\begin{frontmatter}

	\title{A Canonical Internal Model for Disturbance Rejection for a Class of Nonlinear Systems Subject to Trigonometric-Polynomial Disturbances\thanksref{footnoteinfo}} % Title, preferably not more

	\thanks[footnoteinfo]{
		% The work described in this paper was  substantially  supported by a
		% grant from the Research Grants Council of the Hong Kong Special Administrative
		% Region, China [Project No.: CUHK 14203924].
		A preliminary version of this paper was submitted to IFAC 2026 \citep{He26}.
		Corresponding author: Jie Huang. Tel.: +852-39438473.
		Fax: +852-26036002.}

	\author[First]{Changran He}\ead{hechangran@link.cuhk.edu.hk},
	\author[Second]{Jie Huang}\ead{jhuang@mae.cuhk.edu.hk}

	\address[First]{The School of Electrical and Electronic Engineering, Nanyang Technological University, Singapore, 639798 (e-mail: hechangran@link.cuhk.edu.hk).}
	\address[Second]{Department of Mechanical and Automation Engineering, The Chinese University of Hong Kong, Shatin, N.T., Hong Kong SAR (e-mail: jhuang@mae.cuhk.edu.hk).}

	% keyword list or with the
	% help of the Automatica
	% keyword wizard

	\begin{abstract}     % Abstract of not more than 200 words.
		Trigonometric-polynomial disturbances are among the most commonly encountered disturbances in practice, as they can approximate nearly any periodic signal with unknown period. The most effective method for asymptotically rejecting this class of disturbances is through a dynamic compensator known as an internal model, which transforms the disturbance rejection problem into a stabilization problem for an augmented system. However, existing internal model design approaches rely heavily on the properties of the solution to the regulator equations. An effective internal model can only be constructed when this solution exhibits specific characteristics, such as being polynomial in the exogenous signal.
		For complex nonlinear systems - especially nonautonomous ones - solving the disturbance rejection problem using traditional methods remains challenging.
		In this paper, we propose a novel framework for disturbance rejection in a class of nonautonomous nonlinear systems affected by trigonometric-polynomial disturbances. The core of our approach is the design of a canonical internal model that directly converts the disturbance rejection problem into an adaptive stabilization problem for an augmented system. Unlike conventional methods, this internal model is synthesized directly from the given nonlinear plant and the knowledge of the exosystem, without relying on the solution of the regulator equations. This makes the approach applicable to a significantly broader class of nonautonomous nonlinear systems.
		Furthermore, we develop an adaptive disturbance observer comprising the canonical nonlinear internal model, a Luenberger-type state observer, and a parameter adaptation law. This observer ensures global asymptotic convergence of the disturbance estimate to the true disturbance without requiring persistent excitation (PE). Under the PE condition, both the disturbance estimation error and the parameter estimation error converge exponentially.
		By incorporating the disturbance estimate as a feedforward compensation signal, we establish sufficient conditions for achieving global trajectory tracking and asymptotic disturbance rejection.
		The effectiveness of the proposed method is demonstrated through a numerical simulation of a flexible-joint robotic manipulator.
	\end{abstract}

\end{frontmatter}

\section{Introduction}
Mechanical systems are often subjected to disturbances arising from exogenous environments and unmodeled dynamics, which can cause tracking errors, degraded transient performance, or even instability.
The disturbance rejection problem is a special case of the output regulation problem, which seeks to design a feedback control law such that the closed-loop system's output asymptotically tracks a class of reference signals despite the presence of disturbances - where both references and disturbances are generated by an autonomous exosystem.
Research on the linear output regulation problem began in the 1970s \citep{Davison76, Francis77} and was extended to nonlinear systems in the 1990s \citep{HR90, Isidori90}.
A necessary condition for solvability of the nonlinear problem is the existence of solutions to a set of nonlinear partial differential equations known as the regulator equations \citep{Isidori90}. Decades of research have established the internal model principle as a systematic framework for output regulation \citep{Byrnes97, Huang04, Chen15}. In essence, an internal model is a dynamic compensator that converts the original regulation problem into a stabilization problem for an augmented system consisting of the plant and the internal model \citep{Chen04}.

Successful application of this framework hinges on two key issues: (1) the existence of an internal model for the given plant and exosystem, and (2) the stabilizability of the resulting augmented system. The first issue depends on properties of the regulator equations' solution, while the second has been satisfactorily resolved for various classes of standard autonomous nonlinear systems \citep{Chen15}.

There are two widely used internal models: the canonical linear internal model, effective when the solution to the regulator equations is polynomial in the exosystem state \citep{Nikiforov98}, and a more recent type introduced in \cite{CH14} for rejecting trigonometric-polynomial disturbances. However, neither is applicable when the regulator equations' solution is non-polynomial or when the plant is nonautonomous (in which case the regulator equations may not even be well-defined).

Along this line of research, the global disturbance rejection problem of uncertain nonlinear systems in lower triangular form with a neutrally stable uncertain linear exosystem was studied in \cite{Liu11}.
The problem was converted, via the internal model principle, into an adaptive stabilization problem of an augmented system composed of the original plant and the internal model.
However, the internal model in \cite{Liu11} relies on the solution of the regulator equations associated with the plant and exosystem.
Using the second type of internal model, \cite{CH14} studied attitude tracking and disturbance rejection of uncertain rigid-body systems subject to trigonometric-polynomial disturbances with unknown frequencies.
Unlike the approach in \cite{Liu11}, the internal-model-based approach in \cite{CH14} does not require knowledge of the solution to the regulator equations.
However, the approach in \cite{CH14} is specialized for rigid-body systems and does not directly apply to other nonlinear systems.

In this paper, we address the disturbance rejection problem for a broad class of nonautonomous nonlinear systems subject to trigonometric-polynomial disturbances, for which conventional internal-model-based methods do not apply. To solve this problem, we propose an adaptive disturbance observer that integrates a canonical nonlinear internal model, a Luenberger-type state observer, and a parameter adaptation law.
The core innovation is the canonical nonlinear internal model, which is constructed directly from the plant and the knowledge of the exosystem without solving the regulator equations. This makes it suitable for complex nonlinear and nonautonomous systems. Under mild assumptions, the proposed observer achieves global asymptotic convergence of the disturbance estimate to the true disturbance without requiring persistent excitation (PE). When the PE condition is satisfied, both the disturbance estimation error and the parameter estimation error converge exponentially.
By incorporating the disturbance estimate as a feedforward compensation term, we derive sufficient conditions for global trajectory tracking and asymptotic disturbance rejection.

It is worth noting that the concept of an adaptive disturbance observer was previously employed in \cite{Marino07} for linear time-invariant SISO systems with trigonometric-polynomial exogenous signals - a special case of the systems considered here. However, the approach in \cite{Marino07} does not extend to nonlinear or nonautonomous systems.
Moreover, a major difference between the internal-model-based approach in \cite{CH14} and the approach proposed in this paper is that the disturbance is not explicitly reconstructed in \cite{CH14}, whereas the adaptive disturbance observer proposed here explicitly reconstructs the unknown disturbance.
A preliminary version of this paper, which considers only matched disturbances, is under review \citep{He26}.

The remainder of the paper is organized as follows:
Section \ref{SecPreliminaries} reviews preliminary concepts.
Section \ref{SecObserver} presents the adaptive observer design and analyzes the convergence of estimation errors.
Section \ref{SecDisturbanceRejection} establishes sufficient conditions for global disturbance rejection.
Section \ref{SecFrequency} discusses the recovery of unknown disturbance frequencies.
Section \ref{SecExample} illustrates the approach through a numerical example involving trajectory tracking and disturbance rejection in a flexible-joint robotic manipulator.
Finally, Section \ref{SecConclusion} concludes the paper and discusses future research directions.

\section{Preliminaries}\label{SecPreliminaries}

Consider a nonautonomous nonlinear system described by
\begin{align}\label{system}
	\begin{split}
		\dot{x}_1 & = f_1(t, x, u)         \\
		\dot{x}_2 & = f_2(t, x, u) + d(t),
	\end{split}
\end{align}
where $x_1\in\mathbb{R}^{n_1}$ and $x_2\in\mathbb{R}^{n_2}$ decompose the state variable $x: =\operatorname{col}(x_1, x_2)\in\mathbb{R}^n$ with $n=n_1+n_2$; $u\in\mathbb{R}^{m}$ is the input; $f_1:\mathbb{R}\times\mathbb{R}^n\times\mathbb{R}^{m}\to\mathbb{R}^{n_1}$ and $f_2:\mathbb{R}\times\mathbb{R}^n\times\mathbb{R}^{m}\to\mathbb{R}^{n_2}$ are known nonlinear functions; and $d(t)\in\mathbb{R}^{n_2}$ is the disturbance.
% With a slight abuse of notation, we use the shorthand $f_2(t, x, u) := f_2(t, x)+ u$.

Suppose that the $i$th component of the disturbance is a trigonometric polynomial of the form:
\begin{align}\label{dis_component}
	d_i(t) = a_{i0}+\sum_{j=1}^{\rho_i} a_{ij}\sin(\omega_{ij} t+\phi_{ij}), ~i=1,\ldots,n_2,
\end{align}
where $a_{i0}$ is the unknown step magnitude; $\rho_i$ is a positive integer; $a_{ij}$, $\omega_{ij}$, and $\phi_{ij}$ are the unknown amplitudes, unknown frequencies, and unknown phase angles of the sinusoidal functions, respectively.

Without loss of generality, for each $i=1,\ldots, n_2$, we assume that $\omega_{ij}$, $j=1,\ldots, \rho_i$, are distinct positive numbers.
Hence, each component of $d(t)$ is a combination of a step signal and finitely many sinusoidal signals with distinct frequencies.
Then, for each $i=1,\ldots, n_2$, there exists a positive integer $r_{i}$ such that $d_i(t)$ satisfies the following equation:
\begin{align}\label{characterequ}
	\frac{d^{r_{i}}d_i(t)}{d t^{r_{i}}} & = \beta_{i1} d_i(t)+ \beta_{i2}\frac{d d_i(t)}{d t} + \cdots + \beta_{i r_{i}}\frac{d^{r_{i}-1} d_i(t)}{d t^{r_{i}-1}},
\end{align}
where $\beta_{ij}$, $j=1, \ldots, r_{i}$, are real numbers.

We introduce some auxiliary variables for analytical convenience.
For $i=1, \ldots, n_2$, let
\begin{align}\label{Phimat}
	\begin{split}
		\upsilon_{i}(t) & =	\begin{bmatrix}
			                    d_i(t) & \frac{d d_i(t)}{d t} & \cdots & \frac{d^{r_{i}-1} d_i(t)}{d t^{r_{i}-1}}
		                    \end{bmatrix}^{\top}\in\mathbb{R}^{r_i \times 1} \\
		\Phi_{i}        & =\begin{bmatrix}
			                   0          & 1          & 0          & \cdots & 0             \\
			                   0          & 0          & 1          & \cdots & 0             \\
			                   \vdots     & \vdots     & \vdots     & \ddots & \vdots        \\
			                   0          & 0          & 0          & \cdots & 1             \\
			                   \beta_{i1} & \beta_{i2} & \beta_{i3} & \cdots & \beta_{i r_i}
		                   \end{bmatrix}\in\mathbb{R}^{r_i \times r_i}                                                            \\
		\Gamma_{i}      & =\begin{bmatrix}
			                   1 & 0 & \cdots & 0
		                   \end{bmatrix}\in\mathbb{R}^{1 \times r_i}.
	\end{split}
\end{align}
Then, from \eqref{characterequ}, it can be verified that $\upsilon_{i}(t)$ is governed by
\begin{align}\label{dynupsilon}
	\dot{\upsilon}_{i}(t)= \Phi_{i} \upsilon_{i}(t), ~i=1,\ldots,n_2
\end{align}
and $d_i(t)$ can be generated by
\begin{align}\label{relation_upsilon_d}
	d_i(t)=\Gamma_{i} \upsilon_{i}(t), ~i=1,\ldots,n_2.
\end{align}

Let $(M_{i}, N_{i})$ be any controllable pair of matrices where $M_{i}\in\mathbb{R}^{r_{i}\times r_{i}}$ is Hurwitz and $N_i\in\mathbb{R}^{r_i\times 1}$.
It can be verified that the pair of matrices $(\Phi_{i}, \Gamma_{i})$ is observable, the pair of matrices $(M_{i}, N_{i})$ is controllable, and  $\Phi_{i}$ and $M_{i}$ have no common eigenvalues. Thus, there exists a unique nonsingular matrix $T_{i}\in\mathbb{R}^{r_{i}\times r_{i}}$ satisfying the following Sylvester equation \cite{Nikiforov98}:
\begin{align}\label{sylvesterequ}
	T_{i}\Phi_{i}- M_{i} T_{i} = N_{i} \Gamma_{i}.
\end{align}

For conciseness of notation, define
\begin{align}\label{def_Mmat}
	\begin{split}
		\Phi & :=\operatorname{blkdiag}(\Phi_{1}, \ldots, \Phi_{n_2}), ~
		\Gamma  :=\operatorname{blkdiag}(\Gamma_{1}, \ldots, \Gamma_{n_2}) \\
		M    & :=\operatorname{blkdiag}(M_{1}, \ldots, M_{n_2}), ~
		N       :=\operatorname{blkdiag}(N_{1}, \ldots, N_{n_2})           \\
		T    & :=\operatorname{blkdiag}(T_{1}, \ldots, T_{n_2}), ~
		\upsilon  :=\operatorname{col}(\upsilon_{1}, \ldots, \upsilon_{n_2})
	\end{split}
\end{align}
and
\begin{align}\label{def_varrho}
	\begin{split}
		\varrho & = -T\upsilon, \quad
		\Psi      = \Gamma T^{-1}.
	\end{split}
\end{align}
Then, $\varrho$ is governed by the following uncertain LTI exosystem:
\begin{align}\label{dynvarrho}
	\dot{\varrho}= T\Phi T^{-1}\varrho
\end{align}
and $d(t)$ can be expressed by the output of the exosystem as
\begin{align}\label{d_generate}
	d(t) = -\Psi \varrho(t),
\end{align}
where $-\Psi$ is an unknown output matrix of the exosystem.

We also need to introduce the definition of persistently exciting signals.
\begin{Definition}[Persistently Exciting \cite{Marino95}]\label{DefPE}
	A function $W:[t_0,\infty)\to\mathbb{R}^{w_1\times w_2}$ is said to be persistently exciting (PE) if there exist positive constants $\alpha$ and $T_0$ such that
	\begin{align*}
		\int_t^{t+T_0}W(\tau)W^\top(\tau)d\tau \succeq \alpha I_{w_1}, ~\forall t\ge t_0.
	\end{align*}
\end{Definition}

We rephrase Lemma B.2.3 of \cite{Marino95} as follows, which will be useful in the convergence analysis.
\begin{Lemma}\label{LemPEstability}
	Consider the linear time-varying system
	\begin{subequations}\label{stabilitysys}
		\begin{align}
			\dot{x} & =Ax+\Omega^{\top}(t)z \label{abstractsys1} \\
			\dot{z} & =-\Lambda\Omega(t)Px \label{abstractsys2}
		\end{align}
	\end{subequations}
	where $x\in\mathbb{R}^n$, $z\in\mathbb{R}^p$, $A\in\mathbb{R}^{n\times n}$ is Hurwitz, $P\in\mathbb{R}^{n\times n}$ is a symmetric and positive definite matrix satisfying $A^{\top}P+PA=-Q$ with $Q$ symmetric and positive definite, and $\Lambda\in\mathbb{R}^{p\times p}$ is a symmetric and positive definite matrix.
	If $\|\Omega(t)\|$, $\|\dot{\Omega}(t)\|$ are uniformly bounded and $\Omega(t)\in\mathbb{R}^{p\times n}$ is PE, then the origin of system \eqref{stabilitysys} is exponentially stable.
\end{Lemma}

\begin{Remark}\label{RM_zeroing_polynomial}
	The polynomial $p_{i}(\lambda)=\lambda^{r_{i}}-\beta_{i r_{i}}\lambda^{r_{i}-1}-\cdots-\beta_{i2}\lambda-\beta_{i1}$ is called a \textit{zeroing polynomial} of $d_i(t)$ since \eqref{characterequ} can be rewritten as
	$\frac{d^{r_{i}}d_i(t)}{d t^{r_{i}}}  - \beta_{i1} d_i(t)- \beta_{i2}\frac{d d_i(t)}{d t} - \cdots - \beta_{i r_{i}}\frac{d^{r_{i}-1} d_i(t)}{d t^{r_{i}-1}} =0$.
	A zeroing polynomial of $d_i(t)$ of the least degree is called a \textit{minimal zeroing polynomial} of $d_i(t)$ \citep{Chen15}.
	By Example 2.8 of \cite{Chen15}, suppose $a_{i0}=0$ in \eqref{dis_component}, then $\upsilon_{i}(t)=\operatorname{col}(d_i(t), \frac{d d_i(t)}{d t}, \cdots, \frac{d^{r_{i}-1} d_i(t)}{d t^{r_{i}-1}})$ is PE if $r_i\le 2\rho_i$.
	Noting that system \eqref{dynupsilon} is in observable canonical form, if $a_{i0}=0$ and $r_i= 2\rho_i$, then the characteristic polynomial of $\Phi_{i}$ is given by $p_{i}(\lambda):=\lambda^{r_{i}}-\beta_{i r_{i}}\lambda^{r_{i}-1}-\cdots-\beta_{i2}\lambda-\beta_{i1}=\Pi_{j=1}^{\rho_i} (\lambda^2 + \omega_{ij}^2)$ which is also the minimal zeroing polynomial of $d_i(t)$.
	In addition, $\upsilon_{i}(t)$ is PE.
	If $a_{i0}\neq 0$ in \eqref{dis_component} and $r_i= 2\rho_i+1$, then the characteristic polynomial of $\Phi_{i}$ becomes $p_{i}(\lambda):=\lambda^{r_{i}}-\beta_{i r_{i}}\lambda^{r_{i}-1}-\cdots-\beta_{i2}\lambda-\beta_{i1}=\lambda\Pi_{j=1}^{\rho_i} (\lambda^2 + \omega_{ij}^2)$, which is the minimal zeroing polynomial of $d_i(t)$.
	Moreover, $\upsilon_{i}(t)$ is PE.
\end{Remark}

\begin{Remark}\label{RM_PE_implication}
	Consider the trigonometric-polynomial disturbance $d_i(t)$ with distinct frequencies and state $\upsilon_{i}(t)$ defined in \eqref{Phimat}, if $\upsilon_i(t)$ is PE, then system \eqref{dynupsilon} forms a \emph{minimal realization} of the signal $d_i(t)$ in the sense that the state dimension of \eqref{dynupsilon} equals the degree of the minimal zeroing polynomial of $d_i(t)$.
\end{Remark}

\begin{Remark}\label{RM_upsilon_PE}
	The signal $\upsilon(t)$ is PE if and only if $\varrho(t)$ is PE since $\varrho = -T\upsilon$ with $T$ being nonsingular.
	Recall that $T$ is block-diagonal with $T_i$ being the $i$th diagonal block.
	It follows that $\varrho_i = -T_i \upsilon_i$ with $T_i$ nonsingular, and hence $\varrho_i(t)$ is PE if and only if $\upsilon_i(t)$ is PE.
	However, one should note that $\varrho_i(t)$, $i=1,\ldots,n_2$, being PE individually does not in general imply that the stacked signal $\varrho(t)=\operatorname{col}(\varrho_1(t),\ldots,\varrho_{n_2}(t))$ is PE.
\end{Remark}

\section{An Adaptive Disturbance Observer}\label{SecObserver}

\subsection{A Nonlinear Internal Model}

In this section, we will construct an adaptive observer to estimate the unknown disturbances $d(t)$.
For this purpose, let us first consider a dynamic compensator as follows:
\begin{align}\label{internalmodel}
	\dot{\eta} = M \eta + N f_2(t, x, u) - MNx_2,
\end{align}
where $\eta\in\mathbb{R}^{\bar{r}}$ with $\bar{r}=\sum_{i=1}^{n_2}r_{i}$.

As shown in the following lemma, the solution of \eqref{internalmodel} can be used to produce an exponentially convergent estimate of the state $\varrho$ of the uncertain LTI exosystem \eqref{dynvarrho}.
\begin{Lemma}\label{Lem_varrhohat_converge}
	Consider system \eqref{system} subject to external disturbance $d(t)$ whose components are of the form \eqref{dis_component} and dynamic compensator \eqref{internalmodel}.
	Let
	\begin{align}\label{def_varrhohat}
		\hat{\varrho} & := \eta- Nx_2.
	\end{align}
	Then, for any initial conditions $x(0)$ and $\eta(0)$, the solution of \eqref{internalmodel} satisfies
	$\lim_{t\to \infty}\left(\hat{\varrho}(t)- \varrho(t)\right)  = 0$
	exponentially, where $\varrho(t)$ is defined in \eqref{def_varrho}.
\end{Lemma}

\begin{Proof}
	Let $\bar{\eta} := \eta - \varrho$.
	Noting \eqref{sylvesterequ}, from \eqref{dynvarrho} and \eqref{internalmodel}, we have
	\begin{align*}
		\dot{\bar{\eta}} & = M \eta + N f_2(t, x, u) - MNx_2 - \dot{\varrho} \notag                                   \\
		                 & = M \eta + N f_2(t, x, u) - MNx_2 - T\Phi T^{-1}\varrho \notag                             \\
		                 & = M \eta + N f_2(t, x, u) - MNx_2 - \left(M + N\Psi\right)\varrho \notag                   \\
		                 & = M (\bar{\eta} + \varrho) + N f_2(t, x, u) - MNx_2 - \left(M + N\Psi\right)\varrho \notag \\
		                 & = M \bar{\eta} + N f_2(t, x, u) - MNx_2 - N\Psi\varrho.
	\end{align*}
	From \eqref{d_generate}, system \eqref{system} can be rewritten as
	\begin{align}\label{system_rewritten}
		\dot{x}_2 & = f_2(t, x, u) - \Psi\varrho.
	\end{align}

	Let $\tilde{\eta} := \bar{\eta}-Nx_2$.
	Then, we have
	\begin{align}\label{dyn_tilde_eta}
		\dot{\tilde{\eta}} & = M \bar{\eta} + N f_2(t, x, u) - MNx_2 - N\Psi\varrho \notag \\
		                   & \quad - N
		\left(f_2(t, x, u) - \Psi\varrho\right) \notag                                     \\
		                   & = M \bar{\eta} - MNx_2  \notag                                \\
		                   & = M \left(\bar{\eta}-Nx_2\right) \notag                       \\
		                   & = M \tilde{\eta}.
	\end{align}
	Since $M$ is Hurwitz, $\lim_{t\to\infty}\tilde{\eta}(t)=0$ exponentially.
	Since $\hat{\varrho}-\varrho=\tilde{\eta}$, we further have $\lim_{t\to\infty}\left(\hat{\varrho}(t)-\varrho(t)\right)=\lim_{t\to\infty}\tilde{\eta}(t)=0$ exponentially.
\end{Proof}

\begin{Remark}\label{RM_classical_IM_fails}
	As mentioned in the introduction,
	there are two widely used internal models: the canonical linear internal model, which takes the following form \cite{Nikiforov98}:
	\begin{align}\label{classical_IM1}
		\dot{\eta} = M \eta + N u,
	\end{align}
	and the modified canonical linear internal model from \cite{CH14} as follows:
	\begin{align}\label{classical_IM2}
		\dot{\eta} = M \eta + N u - MN x_2.
	\end{align}
	However, neither \eqref{classical_IM1} nor \eqref{classical_IM2} is able to reconstruct the state of the uncertain exosystem \eqref{dynvarrho}.
	Take \eqref{classical_IM2} as an example. Let $
		\hat{\varrho} := \eta - N x_2$ and 	$e := \hat{\varrho} - \varrho$. Then, the error variable $e$ is governed by $\dot{e} = M e - N f_2(t,x,u)$.
	Since the extra term $-N f_2(t,x,u)$ in general does not vanish, the error $e(t)$ here does not converge to zero in general.
	Similarly, we can show that the internal model \eqref{classical_IM1} does not work, either. That is what motivates us to
	conceive \eqref{internalmodel}.
\end{Remark}

\subsection{Disturbance Estimation}
To further estimate the unknown disturbance $d(t)$, let us first consider the special case where the frequencies of $d(t)$ are known.

\begin{Lemma}\label{Lem_d0_converge}
	Consider system \eqref{system}, the external disturbance $d(t)$ whose components are of the form \eqref{dis_component}, and the dynamic compensator \eqref{internalmodel}.
	Suppose, for all $d_i(t)$, $i=1, \ldots, n_2$, the number $\rho_i$ and the frequencies $\omega_{ij}$, $j=1, \ldots, \rho_i$, and whether a constant term is present in \eqref{dis_component} are known.
	Let
	$\hat{d}_0 := -\Psi \hat{\varrho}$, where $\hat{\varrho}$ is defined in \eqref{def_varrhohat}.
	Then, for any initial conditions $x(0)$ and $\eta(0)$, we have
	\begin{align}\label{d0converge}
		\lim_{t\to\infty}\left(\hat{d}_0(t)-d(t)\right)  = 0
	\end{align}
	exponentially.
\end{Lemma}

\begin{Proof}
	Since, for all $i=1,\ldots,n_2$, the number $\rho_i$, the frequencies $\omega_{ij}$,
	$j=1,\ldots,\rho_i$, and whether a constant term is present in \eqref{dis_component}
	(i.e., whether $a_{i0}=0$ or not) are known, the coefficients $\beta_{ij}$, $j=1, \ldots, r_{i}=2 \rho_i+1$, in \eqref{characterequ} can be computed by solving the equation $\lambda^{r_{i}}-\beta_{i r_{i}}\lambda^{r_{i}-1}-\cdots-\beta_{i2}\lambda-\beta_{i1}=\lambda\Pi_{j=1}^{\rho_i} (\lambda^2 + \omega_{ij}^2)$ if $a_{i0}\neq 0$.
	On the other hand, if $a_{i0}=0$, then the coefficients $\beta_{ij}$, $j=1, \ldots, r_{i}=2 \rho_i$, in \eqref{characterequ} can be computed by solving the equation $\lambda^{r_{i}}-\beta_{i r_{i}}\lambda^{r_{i}-1}-\cdots-\beta_{i2}\lambda-\beta_{i1}=\Pi_{j=1}^{\rho_i} (\lambda^2 + \omega_{ij}^2)$.
	Then, the matrix $\Phi_{i}$ in \eqref{Phimat} is known and the matrix $T_{i}$ can be computed directly by solving the Sylvester equation \eqref{sylvesterequ}.
	Since the matrices $T$ and $\Gamma$ are known, we can compute the matrix $\Psi$ using \eqref{def_varrho}.

	By Lemma \ref{Lem_varrhohat_converge}, we have $\lim_{t\to \infty}\left(\hat{\varrho}(t)- \varrho(t)\right)  = 0$ exponentially. Thus, $\lim_{t\to \infty}(\hat{d}_0 (t)- d(t)) = \Psi\lim_{t\to \infty}(-\hat{\varrho}(t)+ \varrho(t)) = 0$ exponentially.
\end{Proof}

\begin{Remark}
	Under the special case where the frequencies of the external disturbance are known, one can effectively estimate the disturbance $d(t)$ using Lemma \ref{Lem_d0_converge}.
	However, the number and the exact values of the frequencies of external disturbances are costly to identify in practice.
	For this reason, the approach in Lemma \ref{Lem_d0_converge} is no longer valid if the number or the exact values of the frequencies are not known precisely.
	Thus, it is desirable to establish a new approach that can estimate the unknown disturbances $d(t)$ without knowing this information.
\end{Remark}

Before introducing the remaining part of the adaptive disturbance observer, let us rewrite $d(t)$ in the following linearly parameterized form for convenience of parameter adaptation:
\begin{align}\label{disturbance_rewrite}
	d(t) & = -\operatorname{col}(\Psi_1 \varrho_1(t), \ldots, \Psi_{n_2} \varrho_{n_2}(t)) \notag \\
	     & = -\begin{bmatrix}
		          \varrho_1^\top(t)\Psi_1^\top \\
		          \vdots                       \\
		          \varrho_{n_2}^\top(t)\Psi_{n_2}^\top
	          \end{bmatrix}\notag                                                 \\
	     & = -\mathcal{M}(\varrho)\mathcal{C}(\Psi),
\end{align}
where $\mathcal{M}(\varrho):=\operatorname{blkdiag}(\varrho_1^\top(t), \ldots, \varrho_{n_2}^\top(t))\in\mathbb{R}^{n_2\times \bar{r}}$ with $\varrho_i = -T_i \upsilon_i$, $i=1,\ldots,n_2$, and $\mathcal{C}(\Psi):=\operatorname{col}(\Psi_1^\top, \ldots, \Psi_{n_2}^\top)\in\mathbb{R}^{\bar{r}}$
with $\Psi_i = \Gamma_i T_i^{-1}$, $i=1,\ldots,n_2$.

Now, we are ready to present the remaining part of the adaptive disturbance observer.
Consider the following dynamic compensator:
\begin{align}
	\dot{\hat{x}}_2 & = f_2(t, x, u)-\mathcal{M}(\hat{\varrho})\theta +K \left(\hat{x}_2-x_2\right)\label{dyn_hatx} \\
	\dot{\theta}    & = \Lambda\mathcal{M}^\top(\hat{\varrho})P\left(\hat{x}_2-x_2\right)\label{dyn_theta},
\end{align}
where $\hat{x}_2\in\mathbb{R}^{n_2}$ is an estimate of the state variable $x_2$;
$\mathcal{M}(\hat{\varrho}):=\operatorname{blkdiag}(\hat{\varrho}_1^\top(t), \ldots, \hat{\varrho}_{n_2}^\top(t))\in\mathbb{R}^{n_2\times \bar{r}}$ is an estimate of $\mathcal{M}(\varrho)$ with $\hat{\varrho}=\operatorname{col}(\hat{\varrho}_1, \ldots, \hat{\varrho}_{n_2})\in\mathbb{R}^{\bar{r}}$ generated by \eqref{internalmodel} via \eqref{def_varrhohat} such that $\hat{\varrho}_i\in\mathbb{R}^{r_i}$, $i=1,\ldots,n_2$; $\theta\in\mathbb{R}^{\bar{r}}$ is the state of the adaptation law; $K\in\mathbb{R}^{n_2\times n_2}$ is Hurwitz; and $\Lambda\in\mathbb{R}^{\bar{r}\times \bar{r}}$ and $P\in\mathbb{R}^{n_2\times n_2}$ are symmetric and positive definite with $P$ satisfying $K^{\top}P+PK=-I_{n_2}$.

In the next theorem, we show that an appropriately designed output of the dynamic compensator composed of \eqref{internalmodel}, \eqref{dyn_hatx}, and \eqref{dyn_theta} can be used to approximate $d(t)$.
\begin{Theorem}\label{Thm_para_converge}
	Consider system \eqref{system}, the external disturbance $d(t)$ whose components are of the form \eqref{dis_component}, and the dynamic compensators \eqref{internalmodel}, \eqref{dyn_hatx}, and \eqref{dyn_theta}.
	Let
	\begin{align}
		\hat{d} := -\mathcal{M}(\hat{\varrho})\theta. \label{def_hatd}
	\end{align}
	Then, for any Hurwitz $K$, any positive definite $\Lambda$, and any initial conditions $x(0)$, $\eta(0)$, $\hat{x}_2(0)$, and $\theta(0)$, we have
	\begin{align}\label{dconverge}
		\lim_{t\to\infty}\left(\hat{x}_2(t)-x_2(t)\right) & = 0,\quad
		\lim_{t\to\infty}\left(\hat{d}(t)-d(t)\right)  = 0.
	\end{align}
	If, in addition, $\varrho_i(t)$ are PE for all $i=1,\ldots,n_2$, then the two equations in \eqref{dconverge} hold exponentially and
	\begin{align}
		\lim_{t\to\infty}\left(\theta(t)-\mathcal{C}(\Psi)\right) & = 0 \label{lim_theta}
	\end{align}
	exponentially.
\end{Theorem}

\begin{Proof}
	Let $\tilde{x}_2 = \hat{x}_2 -x_2$ and $\tilde{\theta} = \theta-\mathcal{C}(\Psi)$.
	From \eqref{system_rewritten}, \eqref{dyn_hatx}, and \eqref{dyn_theta}, the variable $\tilde{x}_2$ is governed by
	\begin{align}\label{dyn_tildex}
		\dot{\tilde{x}}_2 & = \dot{\hat{x}}_2-\dot{x}_2 \notag                                                                                                                                                          \\
		                  & = f_2(t, x, u)-\mathcal{M}(\hat{\varrho})\theta +K \tilde{x}_2 -f_2(t, x, u) + \mathcal{M}(\varrho)\mathcal{C}(\Psi) \notag                                                                 \\
		                  & = K \tilde{x}_2 -\mathcal{M}(\hat{\varrho})\theta   + \mathcal{M}(\varrho)\mathcal{C}(\Psi) \notag                                                                                          \\
		                  & = K \tilde{x}_2 -\mathcal{M}(\hat{\varrho})\theta +\mathcal{M}(\varrho)\mathcal{C}(\Psi)  + \mathcal{M}(\hat{\varrho})\mathcal{C}(\Psi) -\mathcal{M}(\hat{\varrho})\mathcal{C}(\Psi) \notag \\
		                  & = K \tilde{x}_2 -\mathcal{M}(\hat{\varrho})\left(\theta - \mathcal{C}(\Psi) \right)
		-\left(\mathcal{M}(\hat{\varrho})-\mathcal{M}(\varrho)\right)\mathcal{C}(\Psi)\notag                                                                                                                            \\
		                  & = K \tilde{x}_2 -\mathcal{M}(\hat{\varrho})\tilde{\theta}
		-\mathcal{M}(\tilde{\varrho})\mathcal{C}(\Psi),
	\end{align}
	where \eqref{disturbance_rewrite} has been used and $\mathcal{M}(\tilde{\varrho}):=\mathcal{M}(\hat{\varrho})-\mathcal{M}(\varrho)$ with $\tilde{\varrho} := \hat{\varrho}-\varrho$.

	The variable $\tilde{\theta}$ is governed by
	\begin{align}\label{dyn_tildetheta}
		\dot{\tilde{\theta}} & = \dot{\theta}  \notag                                \\
		                     & = \Lambda\mathcal{M}^\top(\hat{\varrho})P\tilde{x}_2.
	\end{align}
	Consider the following auxiliary function:
	\begin{align*}
		V(\tilde{x}_2, \tilde{\theta}) = \tilde{x}_2^{\top}P\tilde{x}_2 + \tilde{\theta}^{\top}\Lambda^{-1}\tilde{\theta}.
	\end{align*}
	Then, the derivative of $V$ along the trajectories of systems \eqref{dyn_tildex} and \eqref{dyn_tildetheta} is
	\begin{align}\label{Vdot}
		\dot{V} & = 2\tilde{x}_2^{\top}P\dot{\tilde{x}}_2 + 2\tilde{\theta}^{\top}\Lambda^{-1}\dot{\tilde{\theta}} \notag \\
		        & = \tilde{x}_2^{\top}\left(K^{\top}P+PK \right)\tilde{x}_2 \notag                                        \\
		        & \quad+ 2\tilde{x}_2^{\top}P\left(-\mathcal{M}(\hat{\varrho})\tilde{\theta}
		-\mathcal{M}(\tilde{\varrho})\mathcal{C}(\Psi) \right) \notag                                                     \\
		        & \quad + 2\tilde{\theta}^{\top}\left(\mathcal{M}^\top(\hat{\varrho})P\tilde{x}_2\right) \notag           \\
		        & = -\tilde{x}_2^{\top}\tilde{x}_2
		-2\tilde{x}_2^{\top}P\mathcal{M}(\tilde{\varrho})\mathcal{C}(\Psi)  \notag                                        \\
		        & = -\tilde{x}_2^{\top}\tilde{x}_2
		-2\tilde{x}_2^{\top}P\Psi\tilde{\varrho},
	\end{align}
	where the identity $\Psi\tilde{\varrho}=\mathcal{M}(\tilde{\varrho})\mathcal{C}(\Psi)$ has been used.

	Next, we will show that $\tilde{x}_2(t)$ is uniformly bounded.
	In fact, systems \eqref{dyn_tildex} and \eqref{dyn_tildetheta} can be rewritten as
	\begin{align}\label{dyn_LTV}
		\begin{bmatrix}
			\dot{\tilde{x}}_2 \\
			\dot{\tilde{\theta}}
		\end{bmatrix} & =
		A_c(t)
		\begin{bmatrix}
			\tilde{x}_2 \\
			\tilde{\theta}
		\end{bmatrix} +w(t),
	\end{align}
	where
	\begin{align*}
		A_c(t) & :=\begin{bmatrix}
			           K                                       & -\mathcal{M}(\hat{\varrho}) \\
			           \Lambda\mathcal{M}^\top(\hat{\varrho})P & 0
		           \end{bmatrix}, ~
		w(t)    := \begin{bmatrix}
			           -\mathcal{M}(\tilde{\varrho})\mathcal{C}(\Psi) \\
			           0
		           \end{bmatrix}.
	\end{align*}
	It can be seen that system \eqref{dyn_LTV} is a linear time-varying system subject to exponentially decaying perturbation $w(t)$.

	Let us first study the stability of the following unperturbed system obtained from \eqref{dyn_LTV} by setting $w(t)=0$:
	\begin{align}\label{dyn_LTV_unperturbed}
		\begin{bmatrix}
			\dot{\tilde{x}}_2 \\
			\dot{\tilde{\theta}}
		\end{bmatrix} & =
		A_c(t)
		\begin{bmatrix}
			\tilde{x}_2 \\
			\tilde{\theta}
		\end{bmatrix}.
	\end{align}
	The derivative of $V(\tilde{x}_2, \tilde{\theta})$ along the trajectories of the unperturbed system \eqref{dyn_LTV_unperturbed} is
	\begin{align*}
		\dot{V} = -\tilde{x}_2^{\top}\tilde{x}_2 \le 0.
	\end{align*}
	Thus, $V(t):=V(\tilde{x}_2(t), \tilde{\theta}(t))$ is uniformly bounded, which implies that the solutions $\tilde{x}_2(t)$ and $\tilde{\theta}(t)$ of the unperturbed system are uniformly bounded.
	Hence, the state transition matrix $\Phi_u(\tau, t)$ of the unperturbed system \eqref{dyn_LTV_unperturbed} is uniformly bounded.

	Next, we consider the perturbed system \eqref{dyn_LTV}.
	Let
	$z(t)  := \operatorname{col}(\tilde{x}_2(t), \tilde{\theta}(t))$.
	Then, the solution of \eqref{dyn_LTV} can be expressed in terms of $\Phi_u(\tau, t)$ as
	\begin{align}\label{solution_z}
		z(t) = \Phi_u(t, 0)z(0)+\int_{0}^{t}\Phi_u(t, \tau)w(\tau)d\tau,
	\end{align}
	where $\Phi_{u}(\tau, t)$ is the state transition matrix of the unperturbed system \eqref{dyn_LTV_unperturbed}.
	Thus, for $t\ge 0$, we have
	\begin{align}\label{zbound}
		\|z(t)\| & \le \|\Phi_u(t, 0)\|\|z(0)\|+\int_{0}^{t}\|\Phi_u(t, \tau)\|\|w(\tau)\|d\tau.
	\end{align}
	Since $\|\Phi_u(t, \tau)\|$ is uniformly bounded and $\lim_{t\to\infty}w(t)=0$ exponentially, the integral $\int_{0}^{t}\|\Phi_u(t, \tau)\|\|w(\tau)\|d\tau$ is uniformly bounded.
	% Thus, $\lim_{t\to\infty}\int_{0}^{t}\|\Phi_u(t, \tau)\|\|w(\tau)\|d\tau$ exists and is finite.
	From \eqref{zbound}, we know that $z(t)$ is also uniformly bounded.
	Thus,
	$\int_{0}^{t}\tilde{x}_2^{\top}(\tau)P\Psi\tilde{\varrho}(\tau)d\tau$
	is uniformly bounded and
	$\lim_{t\to\infty}\int_{0}^{t}\tilde{x}_2^{\top}(\tau)P\Psi\tilde{\varrho}(\tau)d\tau$
	exists and is finite.

	Define another auxiliary function as follows:
	\begin{align}\label{def_W}
		W(t) = V(t) + 2\int_{0}^{t}\tilde{x}_2^{\top}(\tau)P\Psi\tilde{\varrho}(\tau)d\tau.
	\end{align}
	Using \eqref{Vdot}, the derivative of \eqref{def_W} along the trajectories of \eqref{dyn_LTV} is
	\begin{align}\label{Wdot}
		\dot{W} & = \dot{V}+2\tilde{x}_2^{\top}P\Psi\tilde{\varrho} \notag                      \\
		        & = -\tilde{x}_2^{\top}\tilde{x}_2
		-2\tilde{x}_2^{\top}P\Psi\tilde{\varrho}+2\tilde{x}_2^{\top}P\Psi\tilde{\varrho} \notag \\
		        & = -\tilde{x}_2^{\top}\tilde{x}_2 \le 0.
	\end{align}
	Thus, $W(t)\le W(0)$ and is nonincreasing for all $t\ge 0$. Since $W(t)$ is lower bounded, $\lim_{t\to\infty}W(t)$ exists and is finite.

	A direct calculation gives
	\begin{align*}
		\ddot{W} & = -2\tilde{x}_2^{\top}\dot{\tilde{x}}_2 \notag                                      \\
		         & = -2\tilde{x}_2^{\top}\left(K \tilde{x}_2 -\mathcal{M}(\hat{\varrho})\tilde{\theta}
		-\mathcal{M}(\tilde{\varrho})\mathcal{C}(\Psi)  \right)
	\end{align*}
	which is uniformly bounded.
	Thus, $\dot{W}(t)$ is uniformly continuous.
	By Barbalat's lemma, we have $\lim_{t\to\infty}\dot{W}(t)=0$ which implies that $\lim_{t\to\infty}\tilde{x}_2(t)=0$.

	From \eqref{dyn_tildex}, we have
	\begin{align*}
		\ddot{\tilde{x}}_2 & = K \dot{\tilde{x}}_2 -\dot{\mathcal{M}}(\hat{\varrho})\tilde{\theta}-\mathcal{M}(\hat{\varrho})\dot{\tilde{\theta}} -\dot{\mathcal{M}}(\tilde{\varrho})\mathcal{C}(\Psi).
	\end{align*}
	Next, we will show that $\ddot{\tilde{x}}_2(t)$ is uniformly bounded.
	For this purpose, we first show that $\dot{\hat{\varrho}}(t)$ is uniformly bounded.
	In fact, from \eqref{internalmodel}, \eqref{def_varrhohat}, and \eqref{system_rewritten}, we have
	\begin{align}\label{dyn_varrhohat}
		\dot{\hat{\varrho}} & = \dot{\eta}- N\dot{x}_2 \notag                                                      \\
		                    & = M \eta + N f_2(t, x, u) - MNx_2 - N\left( f_2(t, x, u) - \Psi\varrho\right) \notag \\
		                    & = M (\bar{\eta} + \varrho) - MNx_2 - N\left(  - \Psi\varrho\right) \notag            \\
		                    & = M (\bar{\eta}-Nx_2 + \varrho) + N\Psi\varrho \notag                                \\
		                    & = M (\tilde{\eta} + \varrho) + N\Psi\varrho \notag                                   \\
		                    & = M \hat{\varrho} + N\Psi\varrho,
	\end{align}
	where $\bar{\eta} = \eta - \varrho$ and $\tilde{\eta} = \bar{\eta} - Nx_2$.
	Since $\varrho(t)$ is governed by the neutrally stable LTI system \eqref{dynvarrho},
	both $\varrho(t)$ and $\dot{\varrho}(t)$ are uniformly bounded.
	By Lemma \ref{Lem_varrhohat_converge}, $\lim_{t\to \infty}\left(\hat{\varrho}(t)- \varrho(t)\right)  = 0$ exponentially, which implies that $\hat{\varrho}(t)$ is uniformly bounded.
	From \eqref{dyn_varrhohat}, $\dot{\hat{\varrho}}(t)$ is uniformly bounded, which implies $\dot{\tilde{\varrho}}(t)$ and hence $\dot{\mathcal{M}}(\tilde{\varrho}(t))$ are uniformly bounded.
	Thus, $\ddot{\tilde{x}}_2(t)$ is uniformly bounded and $\dot{\tilde{x}}_2(t)$ is uniformly continuous.
	Recalling that $\lim_{t\to\infty}\tilde{x}_2(t)=0$, by Barbalat's lemma, we have $\lim_{t\to\infty}\dot{\tilde{x}}_2(t)=0$.
	Then, from \eqref{dyn_tildex}, we have
	\begin{align*}
		\lim_{t\to\infty}\dot{\tilde{x}}_2(t) & =	\lim_{t\to\infty}\left(K \tilde{x}_2(t) -\mathcal{M}(\hat{\varrho}(t))\theta(t)   + \mathcal{M}(\varrho(t))\mathcal{C}(\Psi) \right) \notag \\
		                                      & = 	\lim_{t\to\infty}\left( -\mathcal{M}(\hat{\varrho}(t))\theta(t)   + \mathcal{M}(\varrho(t))\mathcal{C}(\Psi) \right) \notag                \\
		                                      & = \lim_{t\to\infty}\left( \hat{d}(t) -d(t) \right) \notag                                                                                     \\
		                                      & =0,
	\end{align*}
	which implies the second equation in \eqref{dconverge}.

	In what follows, we will show the second part of the theorem.
	For this purpose, we will apply Lemma \ref{LemPEstability} to show the exponential stability of the unperturbed system \eqref{dyn_LTV_unperturbed}.
	In fact, system \eqref{dyn_LTV_unperturbed} is in the form of \eqref{stabilitysys} with $x=\tilde{x}_2$, $z=\tilde{\theta}$, $A=K$, $\Omega(t)=-\mathcal{M}^\top(\hat{\varrho})$, and $Q=I_{n_2}$.
	Noting that $\Omega(t)\Omega^{\top}(t)=\mathcal{M}^\top(\hat{\varrho})\mathcal{M}(\hat{\varrho})=\operatorname{blkdiag}\left(\hat{\varrho}_1\hat{\varrho}_1^\top,\ldots,\hat{\varrho}_{n_2}\hat{\varrho}_{n_2}^\top\right)$.
	Then, from Definition \ref{DefPE}, the function $\Omega(t)$ is PE if and only if $\hat{\varrho}_i(t)$ are PE for all $i=1,\ldots,n_2$.
	Since, for all $i=1,\ldots,n_2$, $\varrho_i(t)$ are PE by assumption and $\lim_{t\to \infty}\left(\hat{\varrho}_i(t)- \varrho_i(t)\right)  = 0$, $\hat{\varrho}_i(t)$, $i=1,\ldots,n_2$, are PE by Lemma 5 of \cite{He21}.
	Hence, $\Omega(t)$ is PE.
	In addition, $\|\Omega(t)\|$ and $\|\dot{\Omega}(t)\|$ are uniformly bounded since $\hat{\varrho}(t)$ and $\dot{\hat{\varrho}}(t)$ are uniformly bounded. Then, by Lemma \ref{LemPEstability}, the origin of the unperturbed system \eqref{dyn_LTV_unperturbed} is uniformly exponentially stable.

	Now, consider the perturbed system \eqref{dyn_LTV}. Since the perturbation $w(t)$ is exponentially decaying and the origin of the unperturbed system \eqref{dyn_LTV_unperturbed} is uniformly exponentially stable,
	% by Lemma 2.8 of \cite{Cai22}, 
	the solution of \eqref{dyn_LTV} is such that
	$\lim_{t\to \infty}\tilde{x}_2(t)=0$ exponentially and
	$\lim_{t\to \infty}\tilde{\theta}(t)=0$ exponentially.
	This implies that $\lim_{t\to\infty}\left(\hat{d}(t)-d(t)\right)  = 0$ exponentially.
\end{Proof}

\begin{Remark}
	The overall adaptive observer composed of \eqref{internalmodel}, \eqref{dyn_hatx}, and \eqref{dyn_theta} can be put together as follows:
	\begin{align}\label{observer_full}
		\begin{split}
			\dot{\eta}      & = M \eta + N f_2(t, x, u) - MNx_2                                             \\
			\dot{\hat{x}}_2 & = f_2(t, x, u)-\mathcal{M}(\hat{\varrho})\theta +K \left(\hat{x}_2-x_2\right) \\
			\dot{\theta}    & = \Lambda\mathcal{M}^\top(\hat{\varrho})P\left(\hat{x}_2-x_2\right),
		\end{split}
	\end{align}
	where $\hat{\varrho}  = \eta- Nx_2$.
	By Theorem \ref{Thm_para_converge}, $\hat{d} = -\mathcal{M}(\hat{\varrho})\theta$ is an estimate of the external disturbance $d$.
	Thus, system \eqref{observer_full} can be viewed as an adaptive observer for the external disturbance $d$.
	In addition, system \eqref{observer_full} is nonlinear in its state variables $\eta$, $\hat{x}_2$, and $\theta$ since $\hat{\varrho}  = \eta- Nx_2$.
	Moreover, the adaptive observer \eqref{observer_full} satisfies a separation principle in the sense that its design is independent of the control input $u$.
\end{Remark}

\begin{Remark}
	By Theorem \ref{Thm_para_converge}, we have $\lim_{t\to\infty}(\hat{d}(t)-d(t))=0$ regardless of whether $\varrho_i(t)$, $i=1,\ldots,n_2$, are PE.
	In other words, the proposed adaptive observer can estimate the unknown disturbances $d(t)$ even if none of $\varrho_i(t)$, $i=1,\ldots,n_2$, are PE.
	However, if $\varrho_i(t)$, $i=1,\ldots,n_2$, are PE, exponential convergence of $\hat{d}(t)$ to $d(t)$ and exponential convergence of $\theta(t)$ to $\mathcal{C}(\Psi)$ can be further guaranteed.
\end{Remark}

\begin{Remark}\label{RM_dimension}
	The adaptive observer \eqref{observer_full} does not require prior knowledge of the exact number $\rho_i$ of distinct sinusoidal modes in each disturbance component $d_i$. One may select $r_i$ (the dimension of the $i$th block $\eta_i$ of $\eta$ in \eqref{internalmodel}) sufficiently large such that $r_i \ge 2\rho_i+1$, which is chosen not less than the degree of the minimal zeroing polynomial of $d_i(t)$ (cf. Remark \ref{RM_zeroing_polynomial}). Nevertheless, choosing $r_i$ strictly larger than this minimal value induces overmodeling of $d_i(t)$ \cite{Marino11}, which increases the order of the adaptive observer and in turn may destroy the PE property of $\hat{\varrho}_i$ and thus prevent exponential (while retaining asymptotic) convergence. Thus, excessively conservative choices of $r_i$ should be avoided.
	An approach to identify the number of the distinct frequencies in $d(t)$ was presented in Lemma 3.1 of \cite{Marino07}.
\end{Remark}

\subsection{Estimation of Derivatives of the Disturbance}

In this subsection, we show that the proposed adaptive disturbance observer can be further used to estimate the derivatives of the disturbance $d(t)$.
Estimates for derivatives of the disturbance are often required in disturbance rejection for unmatched disturbances as will be shown shortly in Section \ref{SecDisturbanceRejection}.

The estimates of the derivatives of $d(t)$ can be defined recursively as follows.
Specifically, let
\begin{align}
	\varepsilon_1 := M \hat{\varrho} - N\hat{d}\in\mathbb{R}^{\bar{r}}.\label{def_varepsilon_1}
\end{align}
For $k=1,2,\ldots$, let
\begin{align}
	\varepsilon_{k+1} & := M \varepsilon_{k} - N \delta_{k}\in\mathbb{R}^{\bar{r}},\label{def_varepsilon_i}
\end{align}
where, for $k=1,2,\ldots$,
\begin{align}\label{def_delta_i}
	\delta_k & := -\mathcal{M}(\varepsilon_{k})\theta\in\mathbb{R}^{n_2},
\end{align}
in which $\mathcal{M}(\varepsilon_{k}):=\operatorname{blkdiag}(\varepsilon_{k,1}^\top, \ldots, \varepsilon_{k,n_2}^\top)\in\mathbb{R}^{n_2\times \bar{r}}$ and $\varepsilon_{k,i}\in\mathbb{R}^{r_i}$, $i=1,\ldots,n_2$, is obtained from the partition $\varepsilon_{k} := \operatorname{col}(\varepsilon_{k,1}, \ldots, \varepsilon_{k,n_2})\in\mathbb{R}^{\bar{r}}$.

The next lemma shows that $\varepsilon_k(t)$ is an estimate of $\varrho^{(k)}(t)$ and $\delta_k(t)$ is an estimate of $d^{(k)}(t)$.
\begin{Lemma}\label{Lem_derivative_estimate}
	If $\varrho_i(t)$, $i=1,\ldots,n_2$, are PE, then for any $k=1,2,\ldots$, we have
	\begin{align*}
		\lim_{t\to\infty}(\varepsilon_k(t)-\varrho^{(k)}(t))=0, \\
		\lim_{t\to\infty}(\delta_k(t)-d^{(k)}(t))  = 0
	\end{align*}
	both exponentially.
\end{Lemma}

\begin{Proof}
	We prove the lemma by induction.
	For $k=1$, noting \eqref{sylvesterequ}, from \eqref{dynvarrho}, \eqref{d_generate}, and \eqref{def_varepsilon_1}, we have
	\begin{align}\label{dyn_epsilon_varrho}
		\varepsilon_1 - \dot{\varrho} & = M \hat{\varrho} -N \hat{d} -T\Phi T^{-1}\varrho \notag                                                                         \\
		                              & = M \varrho + M \left(\hat{\varrho}-\varrho\right) - N d - N \left(\hat{d}-d\right) -T\Phi T^{-1}\varrho \notag                  \\
		                              & = \left(M-T\Phi T^{-1}\right)\varrho - N d + M \left(\hat{\varrho}-\varrho\right)  - N \left(\hat{d}-d\right)  \notag            \\
		                              & = -N\Gamma T^{-1}   \varrho - N d + M \left(\hat{\varrho}-\varrho\right)  - N \left(\hat{d}-d\right)  \notag                     \\
		                              & = -N\Gamma T^{-1}   \varrho + N \Gamma T^{-1} \varrho + M \left(\hat{\varrho}-\varrho\right)  - N \left(\hat{d}-d\right)  \notag \\
		                              & =  M \tilde{\varrho} - N \tilde{d},
	\end{align}
	where $\tilde{\varrho} := \hat{\varrho}-\varrho$ and $\tilde{d} := \hat{d}-d$.

	By Lemma \ref{Lem_varrhohat_converge}, we have $\lim_{t\to\infty}\tilde{\varrho}(t)=0$ exponentially.
	Since $\varrho_i(t)$, $i=1,\ldots,n_2$, are PE, by Theorem \ref{Thm_para_converge}, $\lim_{t\to\infty}\tilde{d}(t)=0$ exponentially. Thus,
	\begin{align}\label{lim_varepsilon_1}
		\lim_{t\to\infty} \left(\varepsilon_1(t) - \dot{\varrho}(t)\right)=0
	\end{align}
	exponentially, which implies that $\varepsilon_1(t)$ is uniformly bounded since $\dot{\varrho}(t)$ is uniformly bounded.

	A direct calculation gives
	\begin{align*}
		\delta_1-\dot{d}
		 & =-\mathcal{M}(\varepsilon_{1})\theta+\dot{\mathcal{M}}(\varrho)\mathcal{C}(\Psi)\notag                                                                             \\
		 & =-\mathcal{M}(\varepsilon_{1})\theta+\mathcal{M}(\dot{\varrho})\mathcal{C}(\Psi) \notag                                                                            \\
		 & =-\mathcal{M}(\varepsilon_{1})\theta + \mathcal{M}(\varepsilon_{1})\mathcal{C}(\Psi)\notag                                                                         \\
		 & \quad -\mathcal{M}(\varepsilon_{1})\mathcal{C}(\Psi)+\mathcal{M}(\dot{\varrho})\mathcal{C}(\Psi)  \notag                                                           \\
		 & =-\mathcal{M}(\varepsilon_{1})\left(\theta- \mathcal{C}(\Psi)\right) -\left(\mathcal{M}(\varepsilon_{1})-\mathcal{M}(\dot{\varrho})\right)\mathcal{C}(\Psi) \notag \\
		 & =-\mathcal{M}(\varepsilon_{1})\left(\theta- \mathcal{C}(\Psi)\right) -\mathcal{M}(\varepsilon_{1}-\dot{\varrho})\mathcal{C}(\Psi).
	\end{align*}
	Since $\varrho_i(t)$, $i=1,\ldots,n_2$, are PE, by \eqref{lim_theta}, \eqref{lim_varepsilon_1} and noting that $\mathcal{M}(\varepsilon_{1}(t))$ is uniformly bounded, we have $\lim_{t\to\infty}\left(\delta_1(t)-\dot{d}(t)\right)=0$ exponentially.

	Next, suppose for some $k\ge 1$, it holds that
	\begin{align}
		\lim_{t\to\infty}(\varepsilon_k(t)-\varrho^{(k)}(t))  = 0, \label{lim_varepsilon_i_assume} \\
		\lim_{t\to\infty}(\delta_k(t)-d^{(k)}(t))  = 0 \label{lim_delta_i_assume}
	\end{align}
	both exponentially.
	Then, for $k+1$, from \eqref{dynvarrho} and \eqref{def_varepsilon_i}, we have
	\begin{align}\label{difference_epsilon_varrho}
		\varepsilon_{k+1} - \varrho^{(k+1)} & = M \varepsilon_{k} - N \delta_{k} - T\Phi T^{-1}\varrho^{(k)} \notag                                           \\
		                                    & = M \varrho^{(k)} + M \left(\varepsilon_{k}-\varrho^{(k)}\right) - N d^{(k)} \notag                             \\
		                                    & \quad- N \left(\delta_{k}-d^{(k)}\right) - T\Phi T^{-1}\varrho^{(k)} \notag                                     \\
		                                    & = \left(M - T\Phi T^{-1}\right) \varrho^{(k)} + M \left(\varepsilon_{k}-\varrho^{(k)}\right) - N d^{(k)} \notag \\
		                                    & \quad- N \left(\delta_{k}-d^{(k)}\right)  \notag                                                                \\
		                                    & = -N\Gamma T^{-1} \varrho^{(k)} - N d^{(k)} + M \left(\varepsilon_{k}-\varrho^{(k)}\right) \notag               \\
		                                    & \quad- N \left(\delta_{k}-d^{(k)}\right)  \notag                                                                \\
		                                    & = -N\left(\Gamma T^{-1} \varrho^{(k)} +  d^{(k)}\right) + M \left(\varepsilon_{k}-\varrho^{(k)}\right) \notag   \\
		                                    & \quad- N \left(\delta_{k}-d^{(k)}\right)\notag                                                                  \\
		                                    & =  M \left(\varepsilon_{k}-\varrho^{(k)}\right) - N \left(\delta_{k}-d^{(k)}\right),
	\end{align}
	where the identity $d^{(k)}(t)= -\Gamma T^{-1} \varrho^{(k)}(t)$ obtained from \eqref{d_generate} has been used.
	Then, using \eqref{lim_varepsilon_i_assume} and \eqref{lim_delta_i_assume} in \eqref{difference_epsilon_varrho} gives
	\begin{align}\label{lim_varepsilon_iplus1}
		\lim_{t\to\infty}\left(\varepsilon_{k+1}(t)- \varrho^{(k+1)}(t)\right)=0
	\end{align}
	exponentially.
	This implies that $\varepsilon_{k+1}(t)$ is uniformly bounded since $\varrho^{(k+1)}(t)$ is uniformly bounded.

	In addition, we have
	\begin{align*}
		 & \quad \delta_{k+1}-d^{(k+1)} \notag                                                                                                                                      \\
		 & =-\mathcal{M}(\varepsilon_{k+1})\theta+ \mathcal{M}(\varrho^{(k+1)})\mathcal{C}(\Psi)\notag                                                                              \\
		 & =-\mathcal{M}(\varepsilon_{k+1})\theta + \mathcal{M}(\varepsilon_{k+1})\mathcal{C}(\Psi)\notag                                                                           \\
		 & \quad -\mathcal{M}(\varepsilon_{k+1})\mathcal{C}(\Psi)+\mathcal{M}(\varrho^{(k+1)})\mathcal{C}(\Psi)\notag                                                               \\
		 & =-\mathcal{M}(\varepsilon_{k+1})\left(\theta- \mathcal{C}(\Psi)\right) -\left(\mathcal{M}(\varepsilon_{k+1})-\mathcal{M}(\varrho^{(k+1)})\right)\mathcal{C}(\Psi) \notag \\
		 & =-\mathcal{M}(\varepsilon_{k+1})\left(\theta- \mathcal{C}(\Psi)\right) -\mathcal{M}(\varepsilon_{k+1}-\varrho^{(k+1)})\mathcal{C}(\Psi).
	\end{align*}
	By \eqref{lim_theta}, \eqref{lim_varepsilon_iplus1} and noting that $\mathcal{M}(\varepsilon_{k+1}(t))$ is uniformly bounded, we have $\lim_{t\to\infty}\left(\delta_{k+1}(t)-d^{(k+1)}(t)\right)=0$ exponentially.
	The proof is complete by induction.
\end{Proof}

\begin{Remark}
	Lemma \ref{Lem_derivative_estimate} shows that, under the PE condition, the proposed adaptive disturbance observer can estimate not only the disturbance $d(t)$ but also its derivatives $d^{(k)}(t)$, $k = 1,2,\ldots$.
	This property distinguishes the adaptive disturbance observer \eqref{observer_full} from that in \cite{Marino07}, which only estimated the disturbance but not its derivatives.
\end{Remark}

\section{Tracking and Disturbance Rejection}\label{SecDisturbanceRejection}

\subsection{Matched Disturbances}

Since the adaptive observer \eqref{observer_full} can generate an estimate of the external disturbance, the disturbance rejection for system \eqref{system} can be solved as if the external disturbance were known.

Subject to matched disturbances, system \eqref{system} can be put without loss of generality in the following form:
\begin{align}\label{system_matched}
	\begin{split}
		\dot{x}_1 & = f_1(t, x)            \\
		\dot{x}_2 & = f_2(t, x)+ u + d(t),
	\end{split}
\end{align}
where $x_1\in\mathbb{R}^{n_1}$, $x_2\in\mathbb{R}^{n_2}$, and $x =\operatorname{col}(x_1, x_2)\in\mathbb{R}^n$ with $n=n_1+n_2$; $u\in\mathbb{R}^{n_2}$ is the input; $f_1:\mathbb{R}\times\mathbb{R}^n\to\mathbb{R}^{n_1}$ and $f_2:\mathbb{R}\times\mathbb{R}^n\to\mathbb{R}^{n_2}$ are known nonlinear functions; and $d(t)\in\mathbb{R}^{n_2}$ is the matched disturbance with respect to the input $u$.

Suppose the desired trajectory $x_d(t)$ is bounded and sufficiently smooth with bounded derivatives, and define the tracking error as $\tilde{x} := x-x_d$.
We consider the following control law:
\begin{subequations}\label{control_law_matched_full}
	\begin{align}
		\begin{split}
			u          & = \varpi(t, x, \psi, x_d, \dot{x}_d, \ldots, x_d^{{(p)}}) - \hat{d} \\
			\dot{\psi} & = \nu(t, x, \psi, x_d, \dot{x}_d, \ldots, x_d^{{(p)}}, \hat{d})
		\end{split}\label{control_law_matched_1} \\
		\begin{split}
			\dot{\eta}      & = M \eta + N \left(f_2(t, x)+ u\right) - MNx_2                                \\
			\dot{\hat{x}}_2 & = f_2(t, x)+ u-\mathcal{M}(\hat{\varrho})\theta +K \left(\hat{x}_2-x_2\right) \\
			\dot{\theta}    & = \Lambda\mathcal{M}^\top(\hat{\varrho})P\left(\hat{x}_2-x_2\right)
		\end{split}\label{control_law_matched_2}
	\end{align}
\end{subequations}
where $p$ is some positive integer;
$\hat{\varrho} = \eta - N x_2$;
the first part \eqref{control_law_matched_1} is a state-feedback control law with $\psi$ being the state of a possibly required dynamic compensator; and
the second part \eqref{control_law_matched_2} is the adaptive observer of the form \eqref{observer_full} with $\hat{d} = -\mathcal{M}(\hat{\varrho})\theta$ being the estimate of $d$.

The trajectory tracking and disturbance rejection problem for matched disturbances can be formulated as follows.
\begin{Problem}\label{Prob_disturbance_rejection}
	Consider system \eqref{system_matched} and the external disturbance $d(t)$ whose components are of the form \eqref{dis_component}.
	Design a control law of the form \eqref{control_law_matched_full} such that for any initial condition of the closed-loop system, the trajectories of the closed-loop system are bounded and satisfy
	$\lim_{t\to\infty} \tilde{x}(t) = 0$.
\end{Problem}

We then obtain the following sufficient condition for solving Problem \ref{Prob_disturbance_rejection}.

\begin{Theorem}\label{Thm_disturbance_rejection}
	Problem \ref{Prob_disturbance_rejection} is solvable if there exists a control law of the form \eqref{control_law_matched_full} and a continuously differentiable scalar function $V(t,\tilde{x}): [0,\infty) \times \mathbb{R}^n \to \mathbb{R}$ such that, along the trajectories of the closed-loop system composed of \eqref{system_matched} and \eqref{control_law_matched_full}, the function $V(t,\tilde{x})$ satisfies
	\begin{align}
		\underline{\alpha}(\|\tilde{x}\|)  \le V(t,\tilde{x}) \le \bar{\alpha}(\|\tilde{x}\|) \\
		\frac{\partial V}{\partial t} + \frac{\partial V}{\partial \tilde{x}} \left(
		\begin{bmatrix}
			f_1(t, \tilde{x}+x_d) \\
			f_2(t, \tilde{x}+x_d)+\varpi(t) - \tilde{d}
		\end{bmatrix} - \dot{x}_d\right)\notag                                            \\
		\le -\alpha(\|\tilde{x}\|)+\sigma(\|\tilde{d}\|)\label{ISS_inequ_changed}
	\end{align}
	for all $(t, \tilde{x}, \tilde{d}) \in [0, \infty) \times \mathbb{R}^n \times \mathbb{R}^{n_2}$, where $\underline{\alpha}$, $\bar{\alpha}$, and $\alpha$ are class $\mathcal{K}_\infty$ functions;
	$\sigma$ is a class $\mathcal{K}$ function;
	$\varpi(t) :=\varpi(t, x(t), \psi(t), x_d(t), \dot{x}_d(t), \ldots, x_d^{(p)}(t))$;
	and $\tilde{d} := \hat{d}-d$ is the disturbance estimation error.
\end{Theorem}

\begin{Proof}
	The closed-loop system composed of \eqref{system_matched} and the first equation of \eqref{control_law_matched_1} is
	\begin{align}
		\begin{split}
			\dot{x}_1 & = f_1(t, x)                              \\
			\dot{x}_2 & = f_2(t, x)+ \varpi(t)-\hat{d}(t) +d(t).
		\end{split}
	\end{align}
	% where $\varpi(t):=\varpi(t, x(t), \psi(t), x_d(t), \dot{x}_d(t), \ldots, x_d^{(p)}(t))$.
	Since $\tilde{x} = x - x_d$, the closed-loop system can be rewritten in terms of $\tilde{x}$ as
	\begin{align}\label{error_system_matched}
		\dot{\tilde{x}} & =  \begin{bmatrix}
			                     f_1(t, \tilde{x}+x_d) \\
			                     f_2(t, \tilde{x}+x_d)+ \varpi(t) - \tilde{d}
		                     \end{bmatrix}
		- \dot{x}_d.
	\end{align}
	Under the conditions of this theorem, along the trajectories of \eqref{error_system_matched}, the time derivative of $V(t,\tilde{x})$ satisfies
	\begin{align}\label{dotV_ISS_inequ}
		\dot{V}(t,\tilde{x}) \le -\alpha(\|\tilde{x}\|)+\sigma(\|\tilde{d}\|)
	\end{align}
	for all $(t, \tilde{x}, \tilde{d}) \in [0, \infty) \times \mathbb{R}^n \times \mathbb{R}^{n_2}$.

	By Theorem 2.7 in \cite{Chen15}, inequality \eqref{dotV_ISS_inequ} implies that the closed-loop system is input-to-state stable (ISS) with respect to the input $\tilde{d}$ and the state $\tilde{x}$ \citep{Sontag89, Edwards00} and that $V(t, \tilde{x})$ is an ISS-Lyapunov function for the closed-loop system.
	Hence, we have
	\begin{align}\label{state_bound}
		\|\tilde{x}(t)\| & \le \beta(\|\tilde{x}(t_0)\|, t-t_0) + \gamma\left(\sup_{t_0\le \tau \le t}\|\tilde{d}(\tau)\|\right)
	\end{align}
	for all $t\ge t_0$, where $\beta$ is a class $\mathcal{KL}$ function and $\gamma$ is a class $\mathcal{K}$ function.

	Applying Theorem \ref{Thm_para_converge} shows that $\lim_{t\to\infty}\tilde{d}(t)=0$.
	Next, we show that $\lim_{t\to\infty}\tilde{x}(t)= 0$.
	For this purpose, fix any $\epsilon>0$.
	Since the function $\gamma(\cdot)\in\mathcal{K}$ is continuous and strictly increasing, there exists
	$\delta>0$ such that $\gamma(\delta)<\epsilon/2$.
	Since $\lim_{t\to\infty}\tilde{d}(t)=0$, there exists $t_0^\prime\ge t_0$ such that
	\begin{align*}
		\sup_{\tau\ge t_0^\prime}\|\tilde{d}(\tau)\|\le \delta,
	\end{align*}
	which implies that for all $t\ge t_0^\prime$,
	\begin{align*}
		\sup_{t_0^\prime\le\tau\le t}\|\tilde{d}(\tau)\|\le \sup_{\tau\ge t_0^\prime}\|\tilde{d}(\tau)\|\le\delta.
	\end{align*}
	Then, we have
	\begin{align*}
		\gamma\left(\sup_{t_0^\prime\le\tau\le t}\|\tilde{d}(\tau)\|\right)\le \gamma(\delta)<\frac{\epsilon}{2}.
	\end{align*}
	Applying \eqref{state_bound} with $t_0=t_0^\prime$ gives, for all $t\ge t_0^\prime$,
	\begin{align*}
		\|\tilde{x}(t)\|
		 & \le \beta(\|\tilde{x}(t_0^\prime)\|, t-t_0^\prime)+\gamma\left(\sup_{t_0^\prime\le\tau\le t}\|\tilde{d}(\tau)\|\right)\notag \\
		 & < \beta(\|\tilde{x}(t_0^\prime)\|, t-t_0^\prime)+\frac{\epsilon}{2}.
	\end{align*}
	Since $\beta(\cdot,\cdot)\in\mathcal{KL}$, there exists $T^\prime>0$ such that $\beta(\|\tilde{x}(t_0^\prime)\|,T^\prime)<\epsilon/2$.
	Therefore, for all $t\ge t_0^\prime+T^\prime$,
	$\|\tilde{x}(t)\|<\epsilon/2+\epsilon/2=\epsilon$.
	Since $\epsilon>0$ is arbitrary, we have $\lim_{t\to\infty}\tilde{x}(t)= 0$.
\end{Proof}

\subsection{Unmatched Disturbances}

To handle unmatched disturbances, we use not only the estimate $\hat{d}(t)$ of the disturbance $d(t)$ but also the estimates $\delta_k(t)$, $k = 1,\ldots,q$, of its time derivatives up to order $q$, where $q$ is a nonnegative integer.

For $k=0,1,\ldots,q$, let
$\chi_k := \operatorname{col}(d, \dot{d}, \ldots, d^{(k)})$ and
$\hat{\chi}_k := \operatorname{col}(\hat{d}, \delta_1, \ldots, \delta_{k})$ be the estimate of $\chi_k$.
Similar to \eqref{control_law_matched_full}, we consider the following control law for unmatched disturbances:
\begin{subequations}\label{control_law_unmatched_full}
	\begin{align}
		\begin{split}
			u          & = \varpi_u(t, x, \psi, x_{d}, \dot{x}_d,  \ldots, x_d^{(p)}, \hat{\chi}_q) \\
			\dot{\psi} & = \nu_u(t, x, \psi, x_{d}, \dot{x}_d, \ldots, x_d^{(p)}, \hat{\chi}_q)
		\end{split}\label{control_law_unmatched_1} \\
		\begin{split}
			\dot{\eta}      & = M \eta + N f_2(t, x, u) - MNx_2                                              \\
			\dot{\hat{x}}_2 & = f_2(t, x, u) -\mathcal{M}(\hat{\varrho})\theta +K \left(\hat{x}_2-x_2\right) \\
			\dot{\theta}    & = \Lambda\mathcal{M}^\top(\hat{\varrho})P\left(\hat{x}_2-x_2\right)
		\end{split}\label{control_law_unmatched_2}
	\end{align}
\end{subequations}
where $p$ is some positive integer; the first part \eqref{control_law_unmatched_1} is a dynamic state-feedback control law with $\psi$ being the state of a dynamic compensator; the second part \eqref{control_law_unmatched_2} is the adaptive observer \eqref{observer_full} for the external disturbance $d$; and $\hat{\chi}_q = \operatorname{col}(\hat{d}, \delta_1, \ldots, \delta_{q})$ is the estimate of $\chi_q = \operatorname{col}(d, \dot{d}, \ldots, d^{(q)})$ generated by \eqref{control_law_unmatched_2} via \eqref{def_hatd} and \eqref{def_delta_i}.

Similar to the matched disturbance case, we define the tracking error as $\tilde{x} := x-x_d$.
The trajectory tracking and disturbance rejection problem for unmatched disturbances can be stated below.
\begin{Problem}\label{Prob_disturbance_rejection_unmatched}
	Consider system \eqref{system} and the external disturbance $d(t)$ whose components are of the form \eqref{dis_component}.
	Design a control law of the form \eqref{control_law_unmatched_full} such that for any initial condition of the closed-loop system, the trajectories of the closed-loop system are bounded and satisfy
	$\lim_{t\to\infty} \tilde{x}(t) = 0$.
\end{Problem}

\begin{Assumption}\label{ASS_f2rewrite}
	The closed-loop system composed of \eqref{system} and the first equation of \eqref{control_law_unmatched_1} can be written as
	\begin{align}\label{f2rewrite}
		\dot{x}_2 & = f_2(t, x, \varpi_u(t)) + d(t)\notag                                                      \\
		          & = f_2(t, x, \varpi_{u0}(t)) + \Theta(t)\left(\bar{p}(\hat{\chi}_q)-\bar{p}(\chi_q)\right),
	\end{align}
	where
	$\varpi_{u0}(t):=\varpi_u(t, x, \psi, x_{d}, \dot{x}_d,  \ldots, x_d^{(p)}, \chi_q)$ is obtained from $\varpi_u(\cdot)$ by replacing its last argument, $\hat{\chi}_q$, by $\chi_q$;
	$\gamma_1$ is some positive integer;
	$\Theta(t)\in\mathbb{R}^{n_2\times \gamma_1}$ is some regressor matrix function independent of $\chi_q$ and $\hat{\chi}_q$;
	$\bar{p}(\cdot)\in\mathbb{R}^{\gamma_1}$ is some polynomial function whose components are finite-order multivariate polynomials of the components of its argument; and $\varpi_u(t):=\varpi_u(t, x(t), \psi(t), x_{d}(t), \dot{x}_d(t),  \ldots, x_d^{(p)}(t), \hat{\chi}_q(t))$.
\end{Assumption}

\begin{Remark}
	The ideal control $\varpi_{u0}(t)$ is the control input that would be applied if the disturbance and its derivatives up to order $q$ were known.
	Hence, the term $\Theta(t)\left(\bar{p}(\hat{\chi}_q)-\bar{p}(\chi_q)\right)$ in \eqref{f2rewrite} can be viewed as a mismatch term induced by the estimation error between $\hat{\chi}_q$ and $\chi_q$.
	In other words, Assumption \ref{ASS_f2rewrite} states that the mismatch between the closed-loop systems under the ideal control $\varpi_{u0}(t)$ and the certainty-equivalence-based control $\varpi_u(t)$ can be parameterized by the product of the regressor matrix $\Theta(t)$ and the difference of the polynomial functions $\bar{p}(\hat{\chi}_q)$ and $\bar{p}(\chi_q)$.
\end{Remark}

\begin{Remark}
	For system \eqref{system_matched} subject to matched disturbance,
	Assumption \ref{ASS_f2rewrite} holds automatically with
	\begin{align}
		\dot{x}_2 & = f_2(t, x)+ \varpi(t)-\hat{d}(t) + d(t)\notag                                            \\
		          & = f_2(t, x, \varpi_{u0}(t))+ \Theta(t)\left(\bar{p}(\hat{\chi}_0)-\bar{p}(\chi_0)\right),
	\end{align}
	where $f_2(t, x, \varpi_{u0}(t)) :=f_2(t, x)+\varpi_{u0}(t)$ with $\varpi_{u0}(t):=\varpi(t)$ defined in \eqref{control_law_matched_1}, $\gamma_1=n_2$,  $\Theta(t)=-I_{n_2}$, $q=0$, $\chi_0(t)=d(t)$, $\hat{\chi}_0(t)=\hat{d}(t)$, and
	$\bar{p}(\cdot)$ reduces to the identity mapping such that $\bar{p}(\hat{\chi}_0)=\hat{d}$ and $\bar{p}(\chi_0)=d$.
\end{Remark}

We then obtain the following sufficient condition for trajectory tracking and disturbance rejection for unmatched disturbances in general nonlinear systems.

\begin{Theorem}\label{Thm_disturbance_rejection_unmatched}
	Problem \ref{Prob_disturbance_rejection_unmatched} is solvable if there exists a control law of the form \eqref{control_law_unmatched_full} and a continuously differentiable scalar function $V(t,\tilde{x}): [0,\infty) \times \mathbb{R}^n \to \mathbb{R}$ such that Assumption \ref{ASS_f2rewrite} holds and, along the trajectories of the closed-loop system composed of \eqref{system} and \eqref{control_law_unmatched_full}, the function $V(t,\tilde{x})$ satisfies
	\begin{align}
		\underline{\alpha}(\|\tilde{x}\|)  \le V(t,\tilde{x}) \le \bar{\alpha}(\|\tilde{x}\|) \\
		\frac{\partial V}{\partial t} + \frac{\partial V}{\partial \tilde{x}} \left(
		\begin{bmatrix}
			f_1(t, \tilde{x}+x_d, \varpi_u(t)) \\
			f_2(t, \tilde{x}+x_d, \varpi_{u0}(t))+\tilde{\zeta}(t)
		\end{bmatrix} - \dot{x}_d\right)\notag                                 \\
		\le -\alpha(\|\tilde{x}\|)+\sigma(\|\tilde{\zeta}\|)\label{ISS_inequ_changed_unmatched}
	\end{align}
	for all $(t, \tilde{x}, \tilde{\zeta}) \in [0, \infty) \times \mathbb{R}^n \times \mathbb{R}^{n_2}$, where $\underline{\alpha}$, $\bar{\alpha}$, and $\alpha$ are class $\mathcal{K}_\infty$ functions, $\sigma$ is a class $\mathcal{K}$ function, and
	$\tilde{\zeta}(t) := \Theta(t)\left(\bar{p}(\hat{\chi}_q)-\bar{p}(\chi_q)\right)$.
	Suppose the following two conditions hold:
	\begin{enumerate}
		\item $\varrho_i(t)$ are PE for all $i=1,\ldots,n_2$;
		\item there exist $C_\Theta,\alpha_\Theta\ge0$ and $C_e,\alpha_e>0$ with $\alpha_e>\alpha_\Theta$ such that $\|\Theta(t)\|\le C_\Theta e^{\alpha_\Theta t}$ and $\|\bar{p}(\chi_q(t))-\bar{p}(\hat\chi_q(t))\|\le C_e e^{-\alpha_e t}$.
	\end{enumerate}
	Then, for any initial condition of the closed-loop system, the trajectories of the closed-loop system remain bounded and satisfy $\lim_{t\to\infty} \tilde{x}(t)=0$.
\end{Theorem}

\begin{Proof}
	Under Assumption \ref{ASS_f2rewrite}, the closed-loop system composed of \eqref{system} and the first equation of \eqref{control_law_unmatched_1} is
	\begin{align}
		\begin{split}
			\dot{x}_1 & = f_1(t, x, \varpi_u(t))                      \\
			\dot{x}_2 & = f_2(t, x, \varpi_{u0}(t))+\tilde{\zeta}(t).
		\end{split}
	\end{align}
	% where $\varpi(t):=\varpi(t, x(t), \psi(t), x_d(t), \dot{x}_d(t), \ldots, x_d^{(p)}(t))$.
	Since $\tilde{x} = x - x_d$, the closed-loop system can be rewritten in terms of $\tilde{x}$ as
	\begin{align}\label{error_system_unmatched}
		\dot{\tilde{x}} & =  \begin{bmatrix}
			                     f_1(t, \tilde{x}+x_d, \varpi_u(t)) \\
			                     f_2(t, \tilde{x}+x_d, \varpi_{u0}(t))+\tilde{\zeta}(t)
		                     \end{bmatrix}
		- \dot{x}_d.
	\end{align}
	Under the conditions of this theorem, along the trajectories of \eqref{error_system_unmatched}, the time derivative of $V(t,\tilde{x})$ satisfies
	\begin{align}\label{dotV_ISS_inequ_unmatched}
		\dot{V}(t,\tilde{x}) \le -\alpha(\|\tilde{x}\|)+\sigma(\|\tilde{\zeta}\|)
	\end{align}
	for all $(t, \tilde{x}, \tilde{\zeta}) \in [0, \infty) \times \mathbb{R}^n \times \mathbb{R}^{n_2}$.

	By Theorem 2.7 in \cite{Chen15},
	% Then, by Theorem 4.19 of \cite{Khalil02},
	the closed-loop system is ISS with respect to the input $\tilde{\zeta}$ and state $\tilde{x}$.
	Moreover, $V(t, \tilde{x})$ is an ISS-Lyapunov function for the closed-loop system.
	By the definition of ISS, the solution of the closed-loop system satisfies
	\begin{align}\label{state_bound_unmatched}
		\|\tilde{x}(t)\| & \le \beta(\|\tilde{x}(t_0)\|, t-t_0) + \gamma\left(\sup_{t_0\le \tau \le t}\|\tilde{\zeta}(\tau)\|\right)
	\end{align}
	for all $t\ge t_0$, where $\beta$ is a class $\mathcal{KL}$ function and $\gamma$ is a class $\mathcal{K}$ function.
	By Lemma \ref{Lem_derivative_estimate}, if $\varrho_i(t)$ are PE for all $i=1,\ldots,n_2$, then
	$\lim_{t\to\infty}(\hat{\chi}_q(t)-\chi_q(t))  = 0$
	exponentially.
	Then, the second itemized condition further ensures that $\lim_{t\to\infty}\tilde{\zeta}(t)=0$ exponentially.
	The rest of the proof for showing $\lim_{t\to\infty}\tilde{x}(t)=0$ is similar to that of Theorem \ref{Thm_disturbance_rejection} and is thus omitted.
\end{Proof}

\begin{Remark}
	As mentioned in the introduction, the core innovation of our approach is the design of the canonical nonlinear internal model \eqref{internalmodel}, which converts the disturbance rejection problem into an adaptive stabilization problem for an augmented system composed of the original system \eqref{system} and the internal model \eqref{internalmodel}.
	Unlike traditional internal-model-based methods (e.g., \cite{Chen04}), which need to stabilize the augmented system as a whole, our design renders the augmented system lower triangular, such that it can be decomposed into two subsystems that can be stabilized independently. Owing to this lower triangular structure, the adaptive stabilization of the augmented system can be reduced to two subproblems: (1) stabilization of the error system associated with the adaptive disturbance observer \eqref{observer_full}, and (2) input-to-state stabilization of the original system \eqref{system}.
	The first subproblem is fully solved in Theorem \ref{Thm_para_converge}, while the second subproblem can be satisfactorily addressed for several classes of nonlinear systems using existing nonlinear control techniques \citep{Chen15}.
\end{Remark}

\section{Frequency Identification}\label{SecFrequency}

In this section, we present an online frequency identification algorithm for the external disturbance $d(t)$ whose components are of the form \eqref{dis_component}.
Let us first establish the following lemma, which shows how to identify the system matrix $A$ of an uncertain LTI system using estimates of the state and its derivative.
\begin{Lemma}\label{Lem_identify_A}
	Consider an uncertain LTI system:
	\begin{equation}\label{eq:trueLTI}
		\dot x(t)=A x(t)
	\end{equation}
	where $A\in\mathbb{R}^{n\times n}$ is an unknown constant matrix and $x(t)\in\mathbb{R}^n$ is the state.

	For any $T_1>0$, define the following matrices:
	\begin{align*}
		S(t) & :=\int_{t-T_1}^{t} x(\tau)x^\top(\tau) d\tau       \\
		Y(t) & :=\int_{t-T_1}^{t} \dot x(\tau)x^\top(\tau) d\tau.
	\end{align*}
	Then, $Y(t)=AS(t)$ for all $t\ge t_0+T_1$.

	Suppose $x$ and $\dot{x}$ satisfy the following assumptions:
	\begin{enumerate}
		\item[\textup{(A1)}] There exist $\beta>0$ and $T_1>0$ such that
		      $S(t)\succeq \beta I$ for all $t\ge t_0+T_1$. In other words, $x(t)$ is PE.

		\item[\textup{(A2)}]
		      Suppose we have estimates $\hat{x}(t)$ and $\widehat{\dot x}(t)$ of $x(t)$ and $\dot x(t)$, respectively, that converge to the true values exponentially.
		      That is, there exist $c,\alpha>0$ such that, for all $t\ge t_0$,
		      \begin{align*}
			      \|\hat{x}(t)-x(t)\|\le c e^{-\alpha t},\quad
			      \|\widehat{\dot x}(t)-\dot x(t)\|\le c e^{-\alpha t}.
		      \end{align*}

		\item[\textup{(A3)}] There exists a finite constant $\bar{m}$ such that, for all $t\ge t_0+T_1$,
		      \begin{align*}
			      \max\left\{\sup_{\tau\in[t-T_1,t]}\|x(\tau)\|, \sup_{\tau\in[t-T_1,t]}\|\dot x(\tau)\|\right\}\le \bar{m}.
		      \end{align*}
	\end{enumerate}
	Define
	\begin{align*}
		\hat{S}(t) & :=\int_{t-T_1}^{t}\hat{x}(\tau)\hat{x}^\top(\tau) d\tau          \\
		\hat{Y}(t) & :=\int_{t-T_1}^{t}\widehat{\dot x}(\tau)\hat{x}^\top(\tau) d\tau \\
		\hat{A}(t) & :=\hat{Y}(t)\hat{S}^{-1}(t)
	\end{align*}
	whenever $\hat{S}(t)$ is nonsingular.
	Then, there exist constants $k>0$, $\bar\alpha\in(0,\alpha]$, and a time $t_1\ge t_0+T_1$ such that
	\begin{align*}
		\|\hat{A}(t)-A\|\le k e^{-\bar\alpha t}\quad \forall t\ge t_1.
	\end{align*}
\end{Lemma}

\begin{Proof}
	From \eqref{eq:trueLTI}, for any $T_1>0$, it can be verified that
	\begin{align*}
		Y(t)
		 & = \int_{t-T_1}^{t} (A x(\tau))x^\top(\tau) d\tau            \\
		 & = A \left(\int_{t-T_1}^{t} x(\tau)x^\top(\tau) d\tau\right) \\
		 & = AS(t).
	\end{align*}
	Define
	$\tilde{S}(t) := \hat{S}(t)-S(t)$ and
	$\tilde{Y}(t) := \hat{Y}(t)-Y(t)$.
	Let the estimation errors be $\tilde x := \hat{x} - x$ and $\widetilde{\dot x} := \widehat{\dot x}-\dot x$.
	By (A2), for all $\tau \ge t_0$,
	\begin{align*}
		\|\tilde x(\tau)\|\le c e^{-\alpha \tau}, \quad
		\|\widetilde{\dot x}(\tau)\|\le c e^{-\alpha \tau}.
	\end{align*}
	By (A3), for all $\tau\in[t-T_1,t]$,
	\begin{align*}
		\|x(\tau)\| \le \bar{m}, \quad
		\|\dot x(\tau)\| \le \bar{m}.
	\end{align*}
	Moreover, we have
	\begin{align*}
		\|\hat{x}(\tau)\|
		 & \le \|x(\tau)\| + \|\tilde x(\tau)\| \\
		 & \le \bar{m} + c e^{-\alpha \tau}     \\
		 & \le \bar{m} + c := \tilde{m} .
	\end{align*}

	The rest of the proof is divided into three steps.

	\noindent\textit{Step 1: we establish bounds on $\tilde{S}$ and $\tilde{Y}$.}
	Expanding $\tilde{S}$ and using triangle and Cauchy-Schwarz inequalities gives
	\begin{align*}
		\|\tilde{S}(t)\|
		 & = \Big\|\int_{t-T_1}^{t}\left(\tilde x x^\top + x \tilde x^\top + \tilde x \tilde x^\top\right)d\tau\Big\| \\
		 & \le 2\int_{t-T_1}^{t}\|\tilde x(\tau)\|\|x(\tau)\|d\tau
		+ \int_{t-T_1}^{t}\|\tilde x(\tau)\|^2 d\tau                                                                  \\
		 & \le 2\bar{m} \int_{t-T_1}^{t} c e^{-\alpha \tau} d\tau
		+ \int_{t-T_1}^{t} c^{2} e^{-2\alpha \tau} d\tau                                                              \\
		 & \le \frac{2\bar{m} c}{\alpha}\left(1-e^{-\alpha T_1}\right) e^{-\alpha (t-T_1)} \notag                     \\
		 & \quad + \frac{c^{2}}{2\alpha}\left(1-e^{-2\alpha T_1}\right) e^{-2\alpha (t-T_1)}                          \\
		 & \le c_S  e^{-\alpha (t-T_1)}
	\end{align*}
	for some constant $c_S>0$ depending on $(\bar{m},c,\alpha,T_1)$.

	Similarly, for $\tilde{Y}$, we have
	\begin{align*}
		\|\tilde{Y}(t)\|
		 & = \Big\|\int_{t-T_1}^{t}\left(\widetilde{\dot x}\hat{x}^\top + \dot x\tilde x^\top\right)d\tau\Big\|                  \\
		 & \le \int_{t-T_1}^{t}\Big(\|\widetilde{\dot x}(\tau)\|\|\hat{x}(\tau)\| + \|\dot x(\tau)\|\|\tilde x(\tau)\|\Big)d\tau \\
		 & \le \int_{t-T_1}^{t}\Big(c e^{-\alpha \tau}\tilde{m} + \bar{m}  c e^{-\alpha \tau}\Big)d\tau                          \\
		 & = \frac{c(\tilde{m}+\bar{m})}{\alpha}\left(1-e^{-\alpha T_1}\right) e^{-\alpha (t-T_1)}                               \\
		 & \le c_Y  e^{-\alpha (t-T_1)}
	\end{align*}
	for some constant $c_Y>0$.

	\noindent\textit{Step 2: we show positive definiteness of $\hat{S}(t)$ for sufficiently large $t$.}
	By (A1), $S(t)\succeq \beta I$ for all $t\ge t_0+T_1$, hence $\|S^{-1}(t)\|\le \beta^{-1}$.
	Since $\|\tilde{S}(t)\|\to 0$ exponentially, there exists $t_1\ge t_0+T_1$ such that
	\begin{align*}
		\|\tilde{S}(t)\|\le \frac{\beta}{2}, \quad \forall t\ge t_1 .
	\end{align*}
	Then, we have the following three inequalities for all $t\ge t_1$:
	\begin{align*}
		\hat{S}(t)          & = S(t)+\tilde{S}(t)\succeq \frac{\beta}{2}I, \\
		\|\hat{S}^{-1}(t)\| & \le \frac{2}{\beta},                         \\
		\|\hat{S}^{-1}(t)-S^{-1}(t)\|
		                    & = \|S^{-1}(t)\tilde{S}(t)\hat{S}^{-1}(t)\|   \\
		                    & \le \frac{2}{\beta^2}\|\tilde{S}(t)\|,
	\end{align*}
	where the following identity obtained from Section 5.8 of \cite{Horn12} has been used: $\hat{S}^{-1}-S^{-1} = -S^{-1}\tilde{S}\hat{S}^{-1}$.

	\noindent\textit{Step 3: we establish error bound for $\hat{A}(t)$.}
	Using $Y(t)=A S(t)$, we have
	\begin{align*}
		\hat{A}(t)-A
		 & = \hat{Y}(t)\hat{S}^{-1}(t)-Y(t)S^{-1}(t) \\
		 & = \tilde{Y}(t)\hat{S}^{-1}(t)
		+ Y(t)\left(\hat{S}^{-1}(t)-S^{-1}(t)\right).
	\end{align*}
	For the first part, we have
	\begin{align*}
		\|\tilde{Y}(t)\hat{S}^{-1}(t)\|
		\le \|\tilde{Y}(t)\|\|\hat{S}^{-1}(t)\|
		\le \frac{2}{\beta} c_Y  e^{-\alpha (t-T_1)}.
	\end{align*}
	For the second part, since $\|S(t)\|\le \int_{t-T_1}^{t}\|x(\tau)\|^2 d\tau \le T_1 \bar{m}^2$ and $Y(t)=A S(t)$,
	\begin{align*}
		\|Y(t)\left(\hat{S}^{-1}(t)-S^{-1}(t)\right)\|
		 & \le \|A\|\|S(t)\|\|\hat{S}^{-1}(t)-S^{-1}(t)\|                    \\
		 & \le \frac{2\|A\|T_1 \bar{m}^2}{\beta^2} c_S  e^{-\alpha (t-T_1)}.
	\end{align*}
	Combining above two inequalities, for all $t\ge t_1$, we have
	\begin{align*}
		\|\hat{A}(t)-A\|
		 & \le \left(\frac{2c_Y}{\beta}+\frac{2\|A\|T_1 \bar{m}^2}{\beta^2}c_S\right) e^{-\alpha (t-T_1)} \\
		 & := k e^{-\bar\alpha t}
	\end{align*}
	with $\bar\alpha:=\alpha$ and $k>0$ independent of $t$.
\end{Proof}

\begin{Remark}
	Lemma \ref{Lem_identify_A} is inspired by the least-squares parameter identification method presented in Section 1.3 of \cite{Ljung99}. While the classical least-squares approach in \cite{Ljung99} relies on the actual state $x$ and its actual derivative $\dot{x}$, Lemma \ref{Lem_identify_A} instead employs their estimates $\hat{x}$ and $\widehat{\dot{x}}$ to asymptotically identify the system matrix $A$.
\end{Remark}

In what follows, we will apply Lemma \ref{Lem_identify_A} to identify the system matrix $\Phi_i$ of the uncertain exosystem \eqref{dynupsilon}.
% To this end, recall that the $i$th component $d_i(t)$ of the external disturbance can be generated by an uncertain exosystem \eqref{dynupsilon} via \eqref{relation_upsilon_d}.
To this end, for $i=1,\ldots,n_2$, define the estimates of $\upsilon_i(t)$ and $\dot{\upsilon}_i(t)$ as follows:
\begin{align}
	\hat{\upsilon}_i(t) := \begin{bmatrix}
		                       \hat{d}_i(t) & \delta_{1i}(t) & \cdots & \delta_{(r_i-1)i}(t)
	                       \end{bmatrix}^{\top}\in\mathbb{R}^{r_i}, \label{def_hatupsilon} \\
	\widehat{\dot{\upsilon}}_i(t) := \begin{bmatrix}
		                                 \delta_{1i}(t) & \delta_{2i}(t) & \cdots & \delta_{r_i i}(t)
	                                 \end{bmatrix}^{\top}\in\mathbb{R}^{r_i}, \label{def_hatupsilon_dot}
\end{align}
where, for $i=1,\ldots,n_2$, $\hat{d}_i(t)\in\mathbb{R}$ is the $i$th component of $\hat{d}(t)$ generated by the adaptive disturbance observer \eqref{observer_full} via \eqref{def_hatd};
and, for $i=1,\ldots,n_2$, $k=1,\ldots,r_i$, $\delta_{ki}(t)\in\mathbb{R}$ is the $i$th component of $\delta_k(t)$ generated by \eqref{observer_full} via \eqref{def_delta_i}.

\begin{Lemma}\label{Lem_Phimat_estimate}
	Suppose $\varrho_i(t)$, $i=1,\ldots,n_2$, are PE.
	Let
	\begin{align*}
		\hat{S}_i(t) & :=\int_{t-T_1}^{t}\hat{\upsilon}_i(\tau)\hat{\upsilon}_i^\top(\tau) d\tau            \\
		\hat{Y}_i(t) & :=\int_{t-T_1}^{t}\widehat{\dot{\upsilon}}_i(\tau)\hat{\upsilon}_i^\top(\tau) d\tau.
	\end{align*}
	Then, there exists a $T_1>0$ and a time $t_1\ge t_0+T_1$ such that
	$\hat{S}_i(t)$ is nonsingular for all $t\ge t_1$.
	Define
	\begin{align*}
		\hat{\Phi}_i(t) & :=\hat{Y}_i(t)\hat{S}_i^{-1}(t)
	\end{align*}
	for $t\ge t_1$.
	Then, there exist $k>0$ and $\bar\alpha\in(0,\alpha]$ such that
	\begin{align*}
		\|\hat{\Phi}_i(t)-\Phi_i\|\le k e^{-\bar\alpha t}\quad \forall t\ge t_1,
	\end{align*}
	where $\Phi_i$ is the unknown system matrix of \eqref{dynupsilon}.
\end{Lemma}

\begin{Proof}
	Since $\varrho_i(t)$, $i=1,\ldots,n_2$, are PE, by Theorem \ref{Thm_para_converge} and Lemma \ref{Lem_derivative_estimate},
	$\lim_{t\to\infty}\left(\hat{\upsilon}_i(t)-\upsilon_i(t)\right)                  = 0$ and
	$\lim_{t\to\infty}\left(\widehat{\dot{\upsilon}}_i(t)-\dot{\upsilon}_i(t)\right)  = 0$ both
	exponentially.
	By Remark \ref{RM_upsilon_PE}, if $\varrho_i(t)$ is PE, then $\upsilon_i(t)$ is PE.
	Moreover, $\upsilon_i(t)$ and $\dot{\upsilon}_i(t)$ are uniformly bounded by the structure of \eqref{dynupsilon}.
	Hence, Assumptions (A1)-(A3) of Lemma \ref{Lem_identify_A} hold with
	$x(t)=\upsilon_i(t)$, $\hat{x}(t)=\hat{\upsilon}_i(t)$, $\dot{x}(t)=\dot{\upsilon}_i(t)$, and $\widehat{\dot{x}}(t)=\widehat{\dot{\upsilon}}_i(t)$, and the conclusion follows from Lemma \ref{Lem_identify_A}.
\end{Proof}

\begin{Remark}
	With $\varrho_i(t)$, $i=1,\ldots,n_2$, being PE, each system matrix $\Phi_i$ of the uncertain exosystem \eqref{dynupsilon} can be exponentially identified by Lemma \ref{Lem_Phimat_estimate}.
	Consequently, the proposed adaptive disturbance observer \eqref{observer_full} not only estimates the disturbance $d(t)$, but also exponentially reconstructs a minimal realization of the exosystem that generates $d(t)$.
\end{Remark}

Next, we present a lemma on the spectrum of $\Phi_i$ and $-\Phi_i^2$.
\begin{Lemma}\label{Lem_spectrum}
	Let the frequencies $\{\omega_{ij}\}_{j=1}^{\rho_i}$ of $d_i(t)$ defined in \eqref{dis_component} be distinct and strictly positive.
	Suppose $\upsilon_{i}(t)\in\mathbb{R}^{r_i}$ is PE.
	Then, the spectrum of the system matrix $\Phi_i$ of \eqref{dynupsilon} is
	\begin{align*}
		\operatorname{spec}(\Phi_i)=
		\begin{cases}
			\{\pm \mathrm{j}\omega_{i1},\dots,\pm \mathrm{j}\omega_{i\rho_i}\},          & r_i=2\rho_i,   \\[1mm]
			\{0\}\cup\{\pm \mathrm{j}\omega_{i1},\dots,\pm \mathrm{j}\omega_{i\rho_i}\}, & r_i=2\rho_i+1,
		\end{cases}
	\end{align*}
	and hence
	\begin{align}\label{spec_Phi_square}
		\operatorname{spec}(-\Phi_i^2)=
		\begin{cases}
			\{\omega_{i1}^2,\dots,\omega_{i\rho_i}^2\},   & r_i=2\rho_i,   \\[1mm]
			\{0,\omega_{i1}^2,\dots,\omega_{i\rho_i}^2\}, & r_i=2\rho_i+1,
		\end{cases}
	\end{align}
	where each $\omega_{ij}^2$ is an eigenvalue of algebraic and geometric multiplicities $2$, and $0$ has multiplicity $1$ when present.
\end{Lemma}

\begin{Proof}
	Since $\upsilon_{i}(t)\in\mathbb{R}^{r_i}$ is PE, by Remark \ref{RM_zeroing_polynomial}, the characteristic polynomial of $\Phi_i$ is given by $p_{i}(\lambda)=\lambda^{r_{i}}-\beta_{i r_{i}}\lambda^{r_{i}-1}-\cdots-\beta_{i2}\lambda-\beta_{i1}=\Pi_{j=1}^{\rho_i} (\lambda^2 + \omega_{ij}^2)$ with $r_i= 2\rho_i$ for the case where $a_{i0}= 0$.
	When $a_{i0}\neq 0$, the characteristic polynomial of $\Phi_i$ is $p_{i}(\lambda)=\lambda^{r_{i}}-\beta_{i r_{i}}\lambda^{r_{i}-1}-\cdots-\beta_{i2}\lambda-\beta_{i1}=\lambda\Pi_{j=1}^{\rho_i} (\lambda^2 + \omega_{ij}^2)$ with $r_i= 2\rho_i+1$.
	Thus, the eigenvalues of $\Phi_i$ are the distinct points $\{\pm \mathrm{j}\omega_{ij}\}$ and possibly $0$ when $r_i$ is odd.

	Since $\Phi_i$ is diagonalizable over complex space by distinct eigenvalues, each simple pair $\pm \mathrm{j}\omega_{ij}$ yields two independent eigenvectors of $-\Phi_i^2$ at $\omega_{ij}^2$.
	Indeed, if $\Phi_i v_\pm = \pm \mathrm{j}\omega_{ij} v_\pm$ with $v_+,v_-$ linearly independent, then $-\Phi_i^2 v_\pm = \omega_{ij}^2 v_\pm$. Thus, $v_+,v_-$ are also linearly independent eigenvectors of $-\Phi_i^2$ for $\omega_{ij}^2$.
	Hence, the geometric and algebraic multiplicities at $\omega_{ij}^2$ equal $2$.
	Thus, we have \eqref{spec_Phi_square}, where each $\omega_{ij}^2$ has geometric and algebraic multiplicities $2$ and $0$ has multiplicity $1$ when present.
\end{Proof}

The following proposition provides an online algorithm to extract the frequency information from the estimate $\hat{\Phi}_i(t)$ obtained from Lemma \ref{Lem_Phimat_estimate}.

\begin{Proposition}\label{Prop_identify_frequency}
	Suppose $\varrho_i(t)\in\mathbb{R}^{r_i}$, $i=1,\ldots,n_2$, are PE.
	By Lemma \ref{Lem_Phimat_estimate}, the estimate $\hat{\Phi}_i(t)\in\mathbb{R}^{r_i \times r_i}$ satisfies
	\begin{align*}
		\|\hat{\Phi}_i(t)-\Phi_i\|\le k e^{-\alpha t}, \quad t\ge t_1
	\end{align*}
	for some constants $k,\alpha>0$.
	Suppose the frequencies $\{\omega_{ij}\}_{j=1}^{\rho_i}$ of $d_i(t)$ defined in \eqref{dis_component} are distinct and strictly positive.
	Let $B(a,\delta):=\{z\in\mathbb{C}:|z-a|<\delta\}$ and let
	\begin{align*}
		\delta := \frac{1}{2} \min\left\{\min_{p\neq q}|\omega_{ip}^2-\omega_{iq}^2|, \min_j \omega_{ij}^2\right\}>0.
	\end{align*}
	Then, there exists $t_2\ge t_1$ such that, for all $t\ge t_2$, the spectrum of $-\hat{\Phi}_i^2(t)$ splits into disjoint clusters:
	\begin{enumerate}
		\item for each $j=1,\dots,\rho_i$, exactly two eigenvalues $\mu_{j,1}(t),\mu_{j,2}(t)$ lie in the disk $B(\omega_{ij}^2,\delta)$, counting algebraic multiplicity;
		\item if $r_i$ is odd, exactly one eigenvalue lies in the disk $B(0,\delta)$, counting algebraic multiplicity.
	\end{enumerate}
	For each $j=1,\dots,\rho_i$, define
	\begin{align}\label{hat_oemga}
		\hat{\Omega}_{ij}(t):=\frac{\mu_{j,1}(t)+\mu_{j,2}(t)}{2}\in\mathbb{R},
		\quad
		\hat{\omega}_{ij}(t):=\sqrt{\hat{\Omega}_{ij}(t)}.
	\end{align}
	Then, there exist constants $\bar{c}>0$, $\bar{\alpha}\in(0,\alpha]$, and $t_3\ge t_2$ such that, for all $t\ge t_3$,
	\begin{align*}
		|\hat{\omega}_{ij}(t)-\omega_{ij}|\le \bar{c} e^{-\bar{\alpha} t},\quad j=1,\dots,\rho_i.
	\end{align*}
\end{Proposition}

\begin{Proof}
	The norm used in the proof is the induced $2$-norm.
	Since
	\begin{align*}
		\hat{\Phi}_i^2(t)-\Phi_i^2
		 & =(\hat{\Phi}_i(t)-\Phi_i)\hat{\Phi}_i(t)+\Phi_i(\hat{\Phi}_i(t)-\Phi_i),
	\end{align*}
	taking norms and using the triangle inequality gives
	\begin{align*}
		\|\hat{\Phi}_i^2(t)-\Phi_i^2\|
		 & \le \|\hat{\Phi}_i(t)-\Phi_i\|\left(\|\hat{\Phi}_i(t)\|+\|\Phi_i\|\right).
	\end{align*}
	Since $\|\hat{\Phi}_i(t)-\Phi_i\|\le k e^{-\alpha t}$, there exists $M>0$ such that $\sup_{t\ge t_1}\|\hat{\Phi}_i(t)\|\le M$. Hence, there exists $k'>0$ with
	\begin{align*}
		\|\hat{\Phi}_i^2(t)-\Phi_i^2\|\le k' e^{-\alpha t},\quad t\ge t_1,
	\end{align*}
	and therefore
	% \begin{align*}
	$\|-\hat{\Phi}_i^2(t)-(-\Phi_i^2)\|\le k' e^{-\alpha t}$,$t\ge t_1$.
	% \end{align*}

	By Remark \ref{RM_upsilon_PE}, $\upsilon_i(t)$ is PE if and only if $\varrho_i(t)$ is PE.
	By Lemma \ref{Lem_spectrum} and the definition of $\delta$, the open disks $B(\omega_{ij}^2,\delta)$ are pairwise disjoint and, in the odd case, are also disjoint from $B(0,\delta)$.

	By Lemma \ref{Lem_spectrum}, the algebraic and geometric multiplicities of each eigenvalue $\omega_{ij}^2$ of $-\Phi_i^2$ equal $2$ (and $0$ has multiplicity $1$ when present).
	Thus, $-\Phi_i^2$ is diagonalizable over the complex space $\mathbb{C}$, and there exists a nonsingular complex matrix $V$ such that
	$-\Phi_i^2 = V \Lambda V^{-1}$,
	where $\Lambda =\mathrm{diag}\left(\{\omega_{ij}^2\}_{j=1}^{\rho_i} \text{ (each twice)}, 0 \text{ if } r_i \text{ is odd}\right)$.
	Define the error matrix
	$E(t):=-\hat{\Phi}_i^2(t)-(-\Phi_i^2)$.
	By the bound proved above, there exists $k'>0$ with
	$\|E(t)\|\le k' e^{-\alpha t}$, $t\ge t_1$.
	The Bauer-Fike theorem (see, e.g., Theorem 6.3.2 of \cite{Horn12}) states that,
	for every eigenvalue $\tilde\lambda$ of $-\hat{\Phi}_i^2(t)= -\Phi_i^2 + E(t)$, there exists an eigenvalue $\lambda$ of $-\Phi_i^2$ such that
	\begin{align*}
		|\tilde\lambda-\lambda|
		\le \kappa(V)\|E(t)\|,
	\end{align*}
	where $\kappa(V):=\|V\|\|V^{-1}\|$ is the condition number.
	Hence, with $c:=\kappa(V)k'$, we have for all $t\ge t_1$,
	\begin{align*}
		\operatorname{spec}(-\hat{\Phi}_i^2(t))
		\subset
		\bigcup_{\lambda\in\operatorname{spec}(-\Phi_i^2)} B\left(\lambda, c e^{-\alpha t}\right).
	\end{align*}
	In particular, for each $j$ and each eigenvalue $\mu_{j,\ell}(t)$ of $-\hat{\Phi}_i^2(t)$ associated to the cluster near $\omega_{ij}^2$, we have
	\begin{align*}
		\big|\mu_{j,\ell}(t)-\omega_{ij}^2\big| \le c e^{-\alpha t},\qquad \ell=1,2 .
	\end{align*}
	Choose $t_2\ge t_1$ such that $c e^{-\alpha t}<\delta$ for all $t\ge t_2$.
	Then, for all $t\ge t_2$, every eigenvalue of $-\hat{\Phi}_i^2(t)$ lies in exactly one of the pairwise disjoint disks $B(\omega_{ij}^2,\delta)$, $j=1,\dots,\rho_i$, and, if $r_i$ is odd, $B(0,\delta)$.

	Moreover, by the invariance of algebraic multiplicities of eigenvalues under small perturbations (see the discussion around equations (1.16)--(1.19) in Chapter II, Section 4 ``Perturbation of the eigenprojections'' of \cite{Kato95}), the total algebraic multiplicity of the eigenvalues of $-\hat{\Phi}_i^2(t)$ contained in each disk $B(\omega_{ij}^2,\delta)$ equals that of $\omega_{ij}^2$ for all $t\ge t_2$.
	In other words, the total algebraic multiplicity of the eigenvalues inside each $B(\omega_{ij}^2,\delta)$ equals that of $\omega_{ij}^2$, namely $2$.
	Similarly, in the odd case, the total algebraic multiplicity inside $B(0,\delta)$ equals $1$.
	Therefore, for each $j$ there are exactly two eigenvalues $\mu_{j,1}(t),\mu_{j,2}(t)$ of $-\hat{\Phi}_i^2(t)$ in $B(\omega_{ij}^2,\delta)$ counting algebraic multiplicity for all $t\ge t_2$.
	In the odd case, there is exactly one more eigenvalue in $B(0,\delta)$.

	For each $j$ and $t\ge t_2$, define
	% \begin{align*}
	$\hat{\Omega}_{ij}(t):=\frac{1}{2}\left(\mu_{j,1}(t)+\mu_{j,2}(t)\right)$.
	% \end{align*}
	By the previous argument and Bauer-Fike theorem (see, e.g., Theorem 6.3.2 of \cite{Horn12}), for $t\ge t_2$, we have
	\begin{align*}
		|\mu_{j,\ell}(t)-\omega_{ij}^2|\le c e^{-\alpha t},\quad \ell=1,2.
	\end{align*}
	Thus, we have, for all $t\ge t_2$,
	\begin{align*}
		|\hat{\Omega}_{ij}(t)-\omega_{ij}^2|
		 & \le \frac{1}{2}\left(|\mu_{j,1}(t)-\omega_{ij}^2|+|\mu_{j,2}(t)-\omega_{ij}^2|\right)\notag \\
		 & \le c e^{-\alpha t}.
	\end{align*}

	Since $\hat{\Phi}_i(t)$ and $-\hat{\Phi}_i^2(t)$ are real matrices, the two eigenvalues in each cluster are either both real or a complex-conjugate pair.
	Thus, in either case $\hat{\Omega}_{ij}(t)\in\mathbb{R}$ for all $t\ge t_2$. Moreover, because $\omega_{ij}^2>0$ and $\hat{\Omega}_{ij}(t)\to \omega_{ij}^2$, there exists $t_3\ge t_2$ such that
	\begin{align*}
		\hat{\Omega}_{ij}(t)\ge \frac{1}{2}\omega_{ij}^2>0,\quad \forall t\ge t_3.
	\end{align*}

	Define $\hat{\omega}_{ij}(t):=\sqrt{\hat{\Omega}_{ij}(t)}$ for $t\ge t_3$. Then
	\begin{align*}
		|\hat{\omega}_{ij}(t)-\omega_{ij}|
		 & =\frac{|\hat{\Omega}_{ij}(t)-\omega_{ij}^2|}{\hat{\omega}_{ij}(t)+\omega_{ij}}\notag \\
		 & \le \frac{|\hat{\Omega}_{ij}(t)-\omega_{ij}^2|}{\omega_{ij}} \notag                  \\
		 & \le \frac{c}{\omega_{ij}} e^{-\alpha t},
		\quad t\ge t_3.
	\end{align*}
	Let $\omega_{\min}:=\min_{1\le j\le \rho_i}\omega_{ij}>0$, and set $\bar{c}:=c/\omega_{\min}$ and $\bar{\alpha}:=\alpha$. Then, for all $t\ge t_3$, it holds that
	\begin{align*}
		|\hat{\omega}_{ij}(t)-\omega_{ij}|\le \bar{c} e^{-\bar{\alpha} t},\quad j=1,\dots,\rho_i.
	\end{align*}
\end{Proof}

\begin{Remark}
	Under the PE condition on $\varrho_i(t)$, Proposition \ref{Prop_identify_frequency} provides an exponentially convergent online frequency-identification algorithm for estimating the frequencies $\omega_{ij}$ of the disturbance component $d_i(t)$ in \eqref{system}. Our frequency-identification scheme applies to nonlinear nonautonomous systems, whereas most existing approaches (e.g., \cite{Marino07}, \cite{Marino11}) focus on linear time-invariant systems.
\end{Remark}

\section{An Example}\label{SecExample}

In this example, we consider the trajectory tracking and disturbance rejection for a planar robot manipulator with two elastic revolute joints.

Let us first recall the standard form of Euler-Lagrange systems with high-order actuator dynamics as follows:
\begin{subequations}\label{EL_system}
	\begin{align}
		H(q)\ddot{q} + C(q, \dot{q})\dot{q} + G(q) & =\xi_{1}  \label{systemEL} \\
		\begin{split}
			\dot{\xi}_{1}   & =\xi_{2}                     \\
			                & ~ \vdots                     \\
			\dot{\xi}_{m-1} & =\xi_{m}                     \\
			\dot{\xi}_{m}   & = \tau + d \label{systemAct}
		\end{split}
	\end{align}
\end{subequations}
where $q\in \mathbb{R}^{n_2}$ is the generalized position and $\dot{q}\in \mathbb{R}^{n_2}$ is the generalized velocity; $H(q) \in \mathbb{R}^{n_2\times n_2}$ is the symmetric and positive definite inertia matrix; $C(q,\dot{q})\dot{q}\in \mathbb{R}^{n_2}$ denotes the Coriolis and centrifugal force vector; $G(q) \in \mathbb{R}^{n_2}$ is the gravitational force vector; $m$ is a positive integer;
for $i=1,\ldots,m$, $\xi_i\in \mathbb{R}^{n_2}$ is the state of the $i$th integrator; $\tau \in \mathbb{R}^{n_2}$ is the control torque vector; and $d\in \mathbb{R}^{n_2}$ is the external disturbance.

As in \cite{He24}, the equations of motion of the robot manipulator are given by
\begin{subequations}\label{examplesys}
	\begin{align}
		H(q_{1})\ddot{q}_{1} + C(q_{1},\dot{q}_{1})\dot{q}_{1} + g(q_{1}) & =K_{0}(q_{2}-q_{1}) \\
		J\ddot{q}_{2} + K_{0}(q_{2}-q_{1})                                & = u_0 + d_0
	\end{align}
\end{subequations}
which is in the form of \eqref{EL_system} with $m=2$, $q=q_{1}$, $G(q)=g(q_{1})+K_{0} q_{1}$, $\xi_{1}= K_{0} q_{2}$, $\xi_{2}= K_{0} \dot{q}_{2}$, $\tau= K_{0} J^{-1}(u_0-K_{0}(q_{2}-q_{1}))$, and $d=K_{0} J^{-1}d_0$.
System \eqref{examplesys} is then in the form of \eqref{system} with $x_1=\operatorname{col}(q,\dot{q},\xi_1)$, $x_2=\xi_2$, $u=\tau$, $f_1(t,x)=\operatorname{col}(\dot{q},H^{-1}(q)(-C(q,\dot{q})\dot{q}-G(q)+\xi_1),\xi_2)\in\mathbb{R}^{n_1}$ satisfying $n_1 = 3n_2$, and $f_2(t,x)=0\in\mathbb{R}^{n_2}$.

In \eqref{examplesys}, the variables and matrices are defined as follows: $n_2=2$,
$q_{i}=\operatorname{col}(q_{i1}, q_{i2})$, $i=1,2$, $K_{0}=I_{2}$, $J=I_{2}$, and
\begin{align*}
	\begin{split}
		H(q_{1})             & = \begin{bmatrix}
			                         2.35+0.16 \cos (q_{12})   & 0.10 + 0.08 \cos (q_{12}) \\
			                         0.10 + 0.08 \cos (q_{12}) & 0.10                      \\
		                         \end{bmatrix}                            \\
		C(q_{1},\dot{q}_{1}) & = \begin{bmatrix}
			                         -0.08 \sin (q_{12})\dot{q}_{12} & -0.08 \sin (q_{12})(\dot{q}_{11}+\dot{q}_{12}) \\
			                         0.08 \sin (q_{12})\dot{q}_{11}  & 0                                              \\
		                         \end{bmatrix} \\
		g(q_{1})             & =\begin{bmatrix}
			                        38.47 \cos (q_{11})+1.83 \cos (q_{11}+q_{12}) \\
			                        1.83 \cos (q_{11}+q_{12})                     \\
		                        \end{bmatrix},
	\end{split}
\end{align*}
which satisfies $H_{m}I_{2} \preceq H(q_{1}) \preceq H_{M}I_{2}$ for some positive constants $H_{m}$ and $H_{M}$.

Let the desired reference signal be
\begin{align*}
	q_{d}(t)
	=\begin{bmatrix}
		 q_{d1}(t) \\
		 q_{d2}(t)
	 \end{bmatrix}
	=\begin{bmatrix}
		 3\sin (\pi t/100) \\
		 4\cos (2\pi t/100)
	 \end{bmatrix}
\end{align*}
which is bounded and sufficiently smooth with bounded derivatives.

The external disturbance $d_0=\operatorname{col}(d_{01}, d_{02})$ to be rejected is given by
\begin{align*}
	d_{01}(t) & =a_{1}\sin(\omega_{1}t+\phi_{1})+1 \\
	d_{02}(t) & =a_{2}\sin(\omega_{2}t+\phi_{2}),
\end{align*}
where, for $i=1,2$, $a_{i}=0.1\cdot i$, $\omega_{i}=i$, and $\phi_{i}=0.1\cdot i$.
Noting that $d(t)=K_{0} J^{-1}d_0(t)=d_0(t)$, the components of $d(t)$ are of the form \eqref{dis_component}.

\subsection{Disturbance Estimation}

Applying \eqref{observer_full} to system \eqref{EL_system}, we have the adaptive disturbance observer as follows:
\begin{align*}
	% \begin{split}
	\dot{\eta}        & = M \eta+ N\tau  - MN\xi_{2}                                                \\
	\dot{\hat{\xi}}_2 & = \tau-\mathcal{M}(\hat{\varrho})\theta +K \left(\hat{\xi}_2-\xi_{2}\right) \\
	\dot{\theta}      & = \Lambda\mathcal{M}^\top(\hat{\varrho}) P\left(\hat{\xi}_2-\xi_{2}\right),
	% \end{split}
\end{align*}
where the design parameters are $K=-10\cdot I_2$, $P=0.05\cdot I_2$, $\Lambda=500\cdot I_{5}$, and
$M_{1}= \begin{bsmallmatrix}
		0  & 1   & 0  \\
		0  & 0   & 1  \\
		-6 & -11 & -6
	\end{bsmallmatrix}$,
$N_{1}=\begin{bsmallmatrix}
		0 \\
		0 \\
		1
	\end{bsmallmatrix}$,
$M_{2}= \begin{bsmallmatrix}
		0  & 1  \\
		-1 & -2
	\end{bsmallmatrix}$,
$N_{2}=\begin{bsmallmatrix}
		0 \\
		1
	\end{bsmallmatrix}$.

Then, it can be verified that $M_i\in\mathbb{R}^{r_i \times r_i}$ and $N_i\in\mathbb{R}^{r_i \times 1}$ are such that $r_1 = 2\rho_1+1$ and $r_2 = 2\rho_2$, where $\rho_i=1$, $i=1,2$, are the numbers of distinct frequencies of the disturbances.
The dimensions of $M_i$ are equal to the degrees of the minimal zeroing polynomials of $d_i(t)$.
By Remark \ref{RM_zeroing_polynomial}, $\upsilon_i(t)$, $i=1,2$, are PE. Then, by Remark \ref{RM_upsilon_PE}, $\varrho_i(t)$, $i=1,2$, are PE.
Hence, by Theorem \ref{Thm_para_converge},
the estimate $\hat{d} = -\mathcal{M}(\hat{\varrho})\theta$ converges to the actual disturbance $d(t)$ exponentially.
Moreover, by Proposition \ref{Prop_identify_frequency}, the frequencies $\omega_1$ and $\omega_2$ can be exponentially identified by their estimates defined in \eqref{hat_oemga}.

\begin{figure}[!htb]
	\centering
	\includegraphics[width=0.9\hsize]{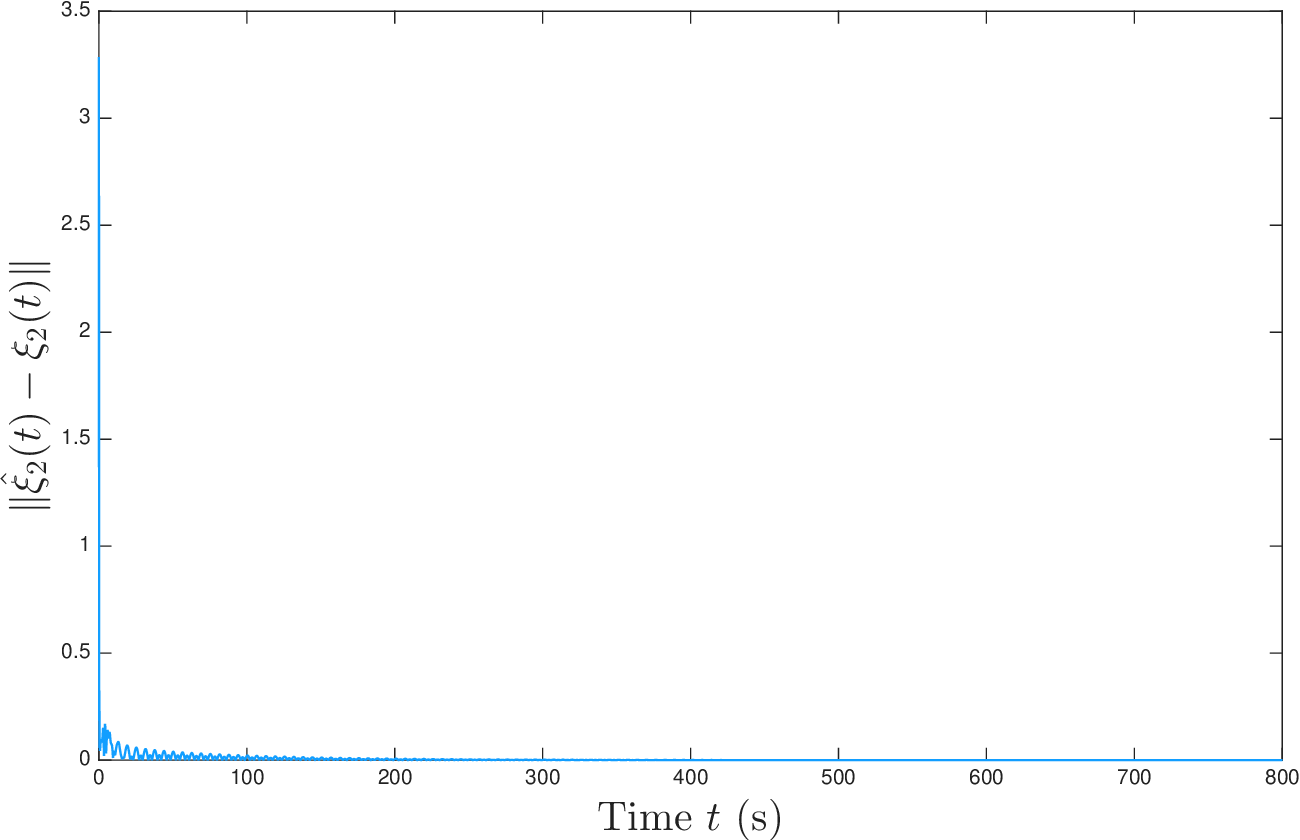}
	\caption{Estimation error of $\hat{\xi}_2(t)$.}
	\label{ehatxi2}
\end{figure}

\begin{figure}[!htb]
	\centering
	\includegraphics[width=0.9\hsize]{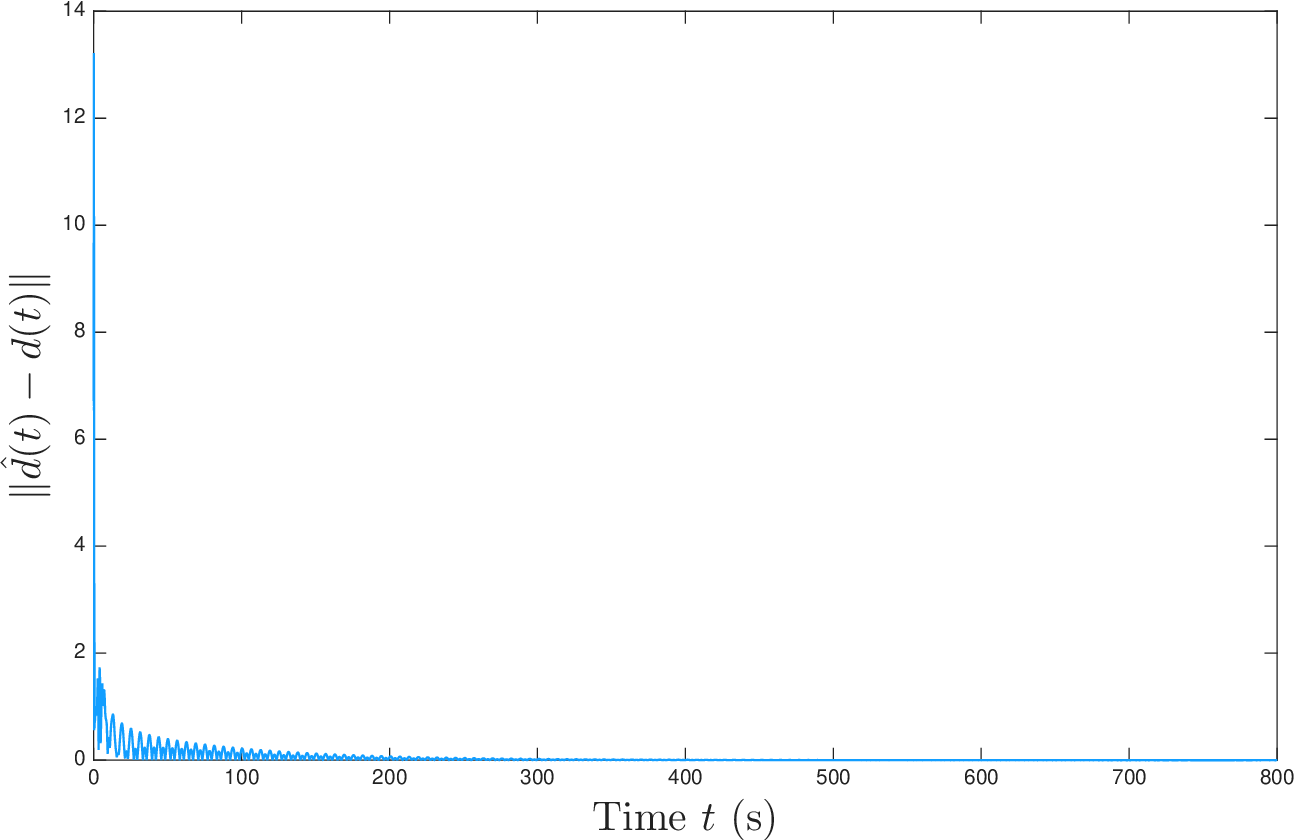}
	\caption{Estimation error of $\hat{d}(t)$.}
	\label{ehatd}
\end{figure}

\begin{figure}[!htb]
	\centering
	\includegraphics[width=0.9\hsize]{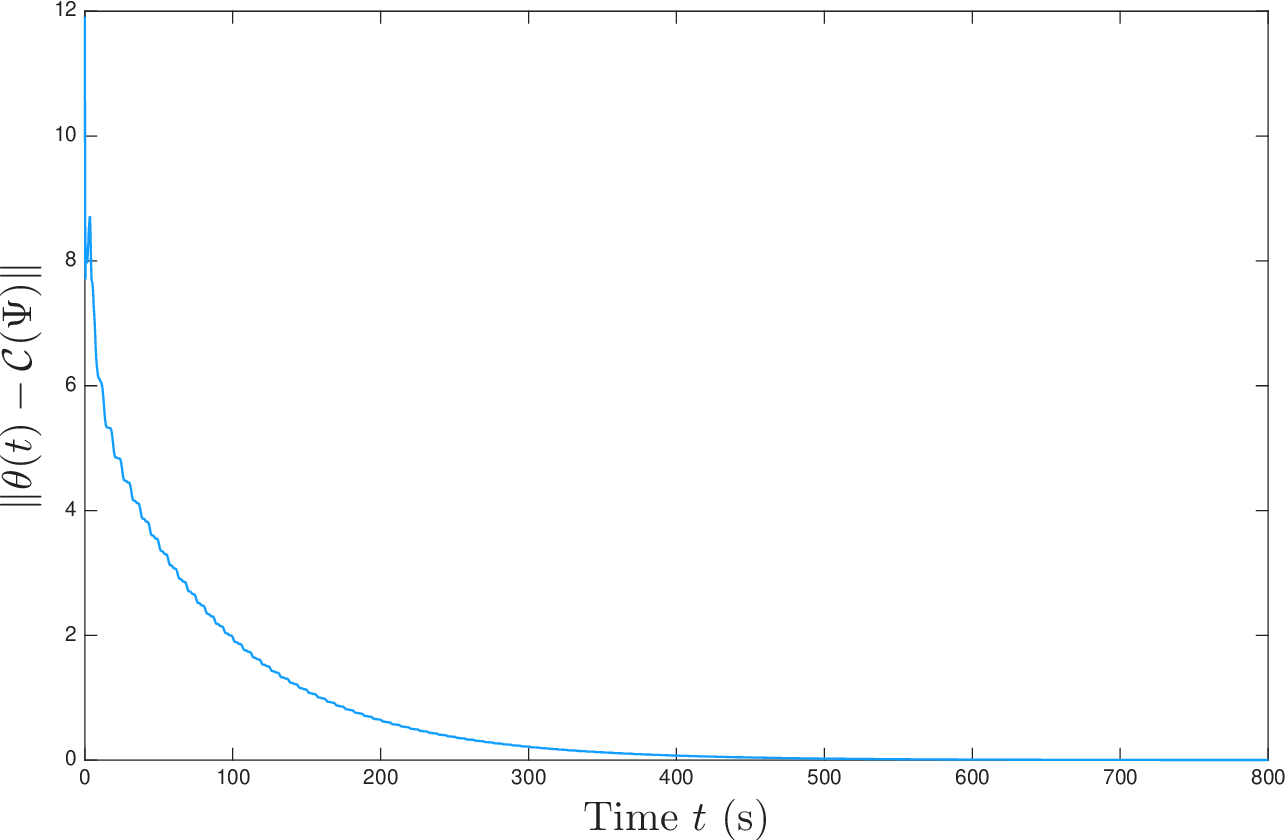}
	\caption{Estimation error of $\theta(t)$.}
	\label{etheta}
\end{figure}

\begin{figure}[!htb]
	\centering
	\includegraphics[width=0.9\hsize]{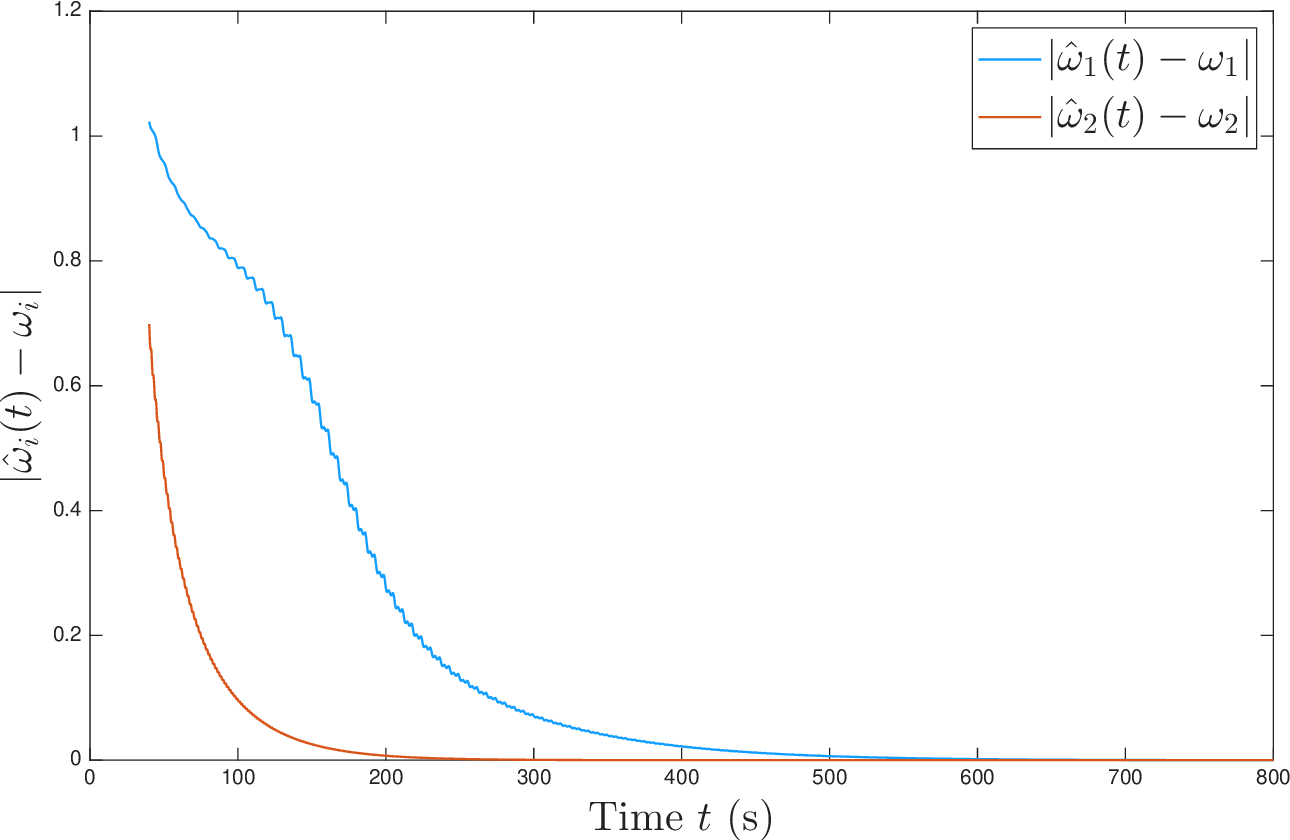}
	\caption{Estimation errors of $\hat{\omega}_1(t)$ and $\hat{\omega}_2(t)$.}
	\label{eomega}
\end{figure}

Computer simulation is conducted with randomly generated initial conditions $q_{i}(0)$, $\dot{q}_{i}(0)$, $\eta(0)$, $\hat{\xi}_2(0)$, $\theta(0)$, and $\xi_{i}(0)$, $i=1,2$, within the interval $[-2, 2]$.
The estimation error of $\hat{\xi}_2(t)$ is shown in Fig. \ref{ehatxi2}, the estimation error of $\hat{d}(t)$ is shown in Fig. \ref{ehatd}, the estimation error of $\theta(t)$ is shown in Fig. \ref{etheta}, and the estimation errors of $\hat{\omega}_1(t)$ and $\hat{\omega}_2(t)$ with $T_1=40$ in Lemma \ref{Lem_Phimat_estimate} are shown in Fig. \ref{eomega}.
Hence, satisfactory disturbance estimation and frequency identification performance are achieved.

\subsection{Disturbance Rejection}

To define the control law, we first define some variables.
Let
$\dot{q}_{r} := \dot{q}_{d} -\alpha (q - q_{d})$
with $\alpha=1$, and let
$s := \dot{q} -\dot{q}_{r}$.
Let
\begin{align}\label{def_xir}
	\xi_r :=-K_s s + y_{0}(q, \dot{q}, \dot{q}_{r}, \ddot{q}_{r}),
\end{align}
where $y_{0}(q, \dot{q}, \dot{q}_{r}, \ddot{q}_{r}):=H(q)\ddot{q}_{r} + C(q, \dot{q})\dot{q}_{r}+ G(q)$
and
$K_s=I_2$.
For $i=1,2$, let $\zeta_{i}  :=\frac{d^{i} \xi_{r}(t)}{d t^{i}}$ be the $i$th derivative of $\xi_{r}$.

Consider the following control law:
\begin{align*}
	\tau & = \zeta_{2} - k_{p1}\left(\xi_1 -\xi_r\right) -k_{p2}\left(\xi_{2}- \zeta_{1}\right) -\hat{d},
\end{align*}
where $k_{p1}=25$ and $k_{p2}=10$.
Let $\tilde{\xi}_{1}:=\xi_1 -\xi_r$, $\tilde{\xi}_{2}  := \xi_{2} -\zeta_{1}$, and $\tilde{\xi}:=\operatorname{col}(\tilde{\xi}_1, \tilde{\xi}_2)$.
By the proof of Theorem 4.1 in \cite{He24}, the closed-loop system is given by
\begin{align}\label{example_closedloop}
	\begin{split}
		H(q) \dot{s}+C(q, \dot{q}) s & =  -K_s s  + C_{c}\tilde{\xi} \\
		\dot{\tilde{\xi}}            & =
		A_{c}\tilde{\xi}-B_{c}\tilde{d},
	\end{split}
\end{align}
where $q=q_1$,
$C_{c} = \begin{bsmallmatrix}
		I_{2} & 0_{2\times 2}
	\end{bsmallmatrix}$,
$A_{c}  =\begin{bsmallmatrix}
		0_{2\times 2} & I_2        \\
		-k_{p1}I_2    & -k_{p2}I_2
	\end{bsmallmatrix}$, and
$B_{c}  =\begin{bsmallmatrix}
		0_{2\times 2} \\
		I_2
	\end{bsmallmatrix}$.
Since $k_{p1}$ and $k_{p2}$ are such that all the roots of $\lambda^{2}+ k_{p1}+ k_{p2}\lambda =0$ have negative real parts, $A_{c}$ is Hurwitz.

\begin{figure}[!htb]
	\centering
	\includegraphics[width=0.9\hsize]{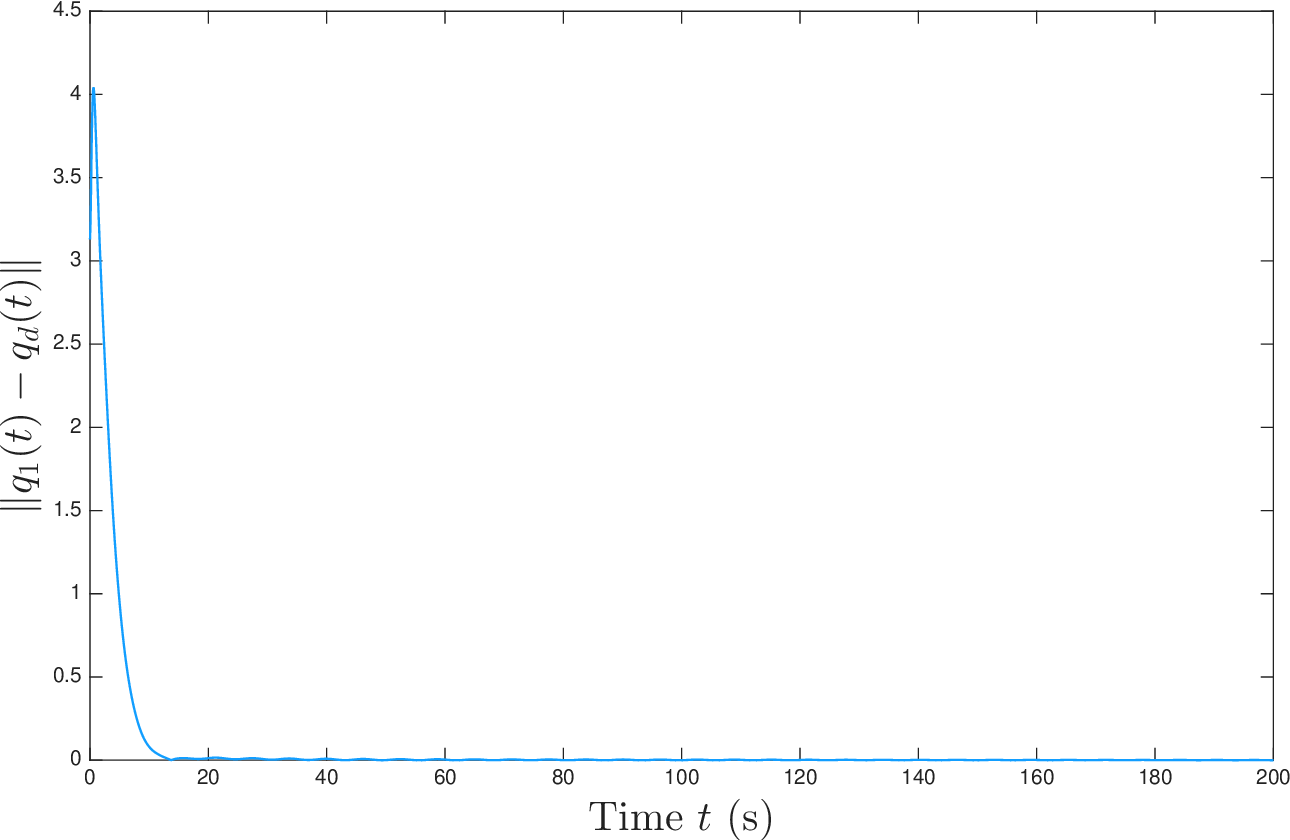}
	\caption{Position tracking error.}
	\label{eq}
\end{figure}

\begin{figure}[!htb]
	\centering
	\includegraphics[width=0.9\hsize]{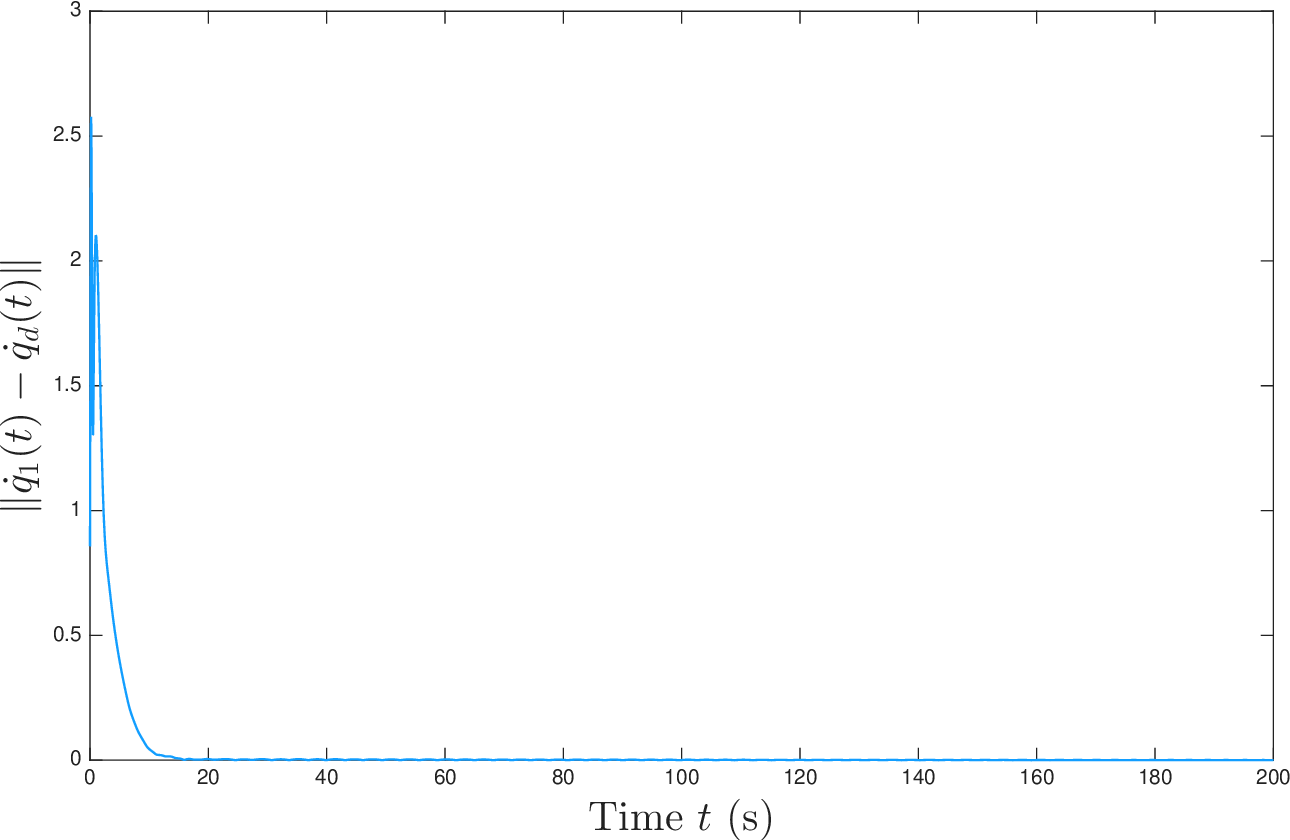}
	\caption{Velocity tracking error.}
	\label{edq}
\end{figure}

Let $P_c\in\mathbb{R}^{4 \times 4}$ be the unique positive definite matrix satisfying $A_{c}^{\top}P_c+ P_c A_{c}=-I_{4}$.
Consider the following Lyapunov function candidate for \eqref{example_closedloop}:
\begin{align*}
	V(t, s, \tilde{\xi}) = \frac{1}{2}s^{\top} H(q(t)) s+c_{1}\tilde{\xi}^{\top}P_c\tilde{\xi},
\end{align*}
where
$c_{1}  >\frac{1}{2\lambda_{\min}(K_s)}+1$.
Then, $V(t, s, \tilde{\xi})$ satisfies
\begin{align*}
	c_{2}\|\operatorname{col}(s, \tilde{\xi})\|^{2} \le V(t, s, \tilde{\xi}) \le c_{3} \|\operatorname{col}(s, \tilde{\xi})\|^{2},
\end{align*}
where $c_{2}:=\min\{\frac{H_{m}}{2}, c_{1}\lambda_{\min}(P_c)\}$, $c_{3}:=\max\{\frac{H_{M}}{2}, c_{1}\lambda_{\max}(P_c)\}$.
Moreover, as shown in the proof of Theorem 4.1 in \cite{He24}, along the trajectories of \eqref{example_closedloop},
\begin{align*}
	\dot{V} & \le   -c_4 V(t, s, \tilde{\xi}) + c_{1}^{2}\|\tilde{d}\|^{2}
\end{align*}
for some constant $c_4>0$.
Thus, all conditions of Theorem \ref{Thm_disturbance_rejection} are satisfied.
By Theorem \ref{Thm_disturbance_rejection}, the state variable $s(t)$ and $\tilde{\xi}(t)$ of the closed-loop system converge to zero asymptotically.
Since $\lim_{t\to\infty}s(t)=0$ implies $\lim_{t\to\infty}(q(t)-q_{d}(t))=0$ and $\lim_{t\to\infty}(\dot{q}(t)-\dot{q}_{d}(t))=0$,
the trajectory tracking and disturbance rejection of system \eqref{examplesys} is achieved.
The position tracking error and the velocity tracking error are shown in Figs. \ref{eq} and \ref{edq}, respectively.
It can be seen that satisfactory tracking performance is achieved.

\begin{Remark}
	Unlike the external disturbance $d(t)$ considered in \cite{He24} whose frequencies are assumed to be known, the frequencies of the components of $d(t)$ in this example are unknown.
\end{Remark}

\begin{Remark}
	It is noteworthy that the proposed control law based on the adaptive disturbance observer achieves trajectory tracking and disturbance rejection without checking the solvability of the regulator equations for the nonlinear plant, which is necessary for solving the conventional nonlinear output regulation problems.
	In general, regulator equations are partial differential equations that are difficult to solve analytically.
\end{Remark}

\section{Conclusion}\label{SecConclusion}

In this paper, we studied disturbance rejection for general nonautonomous nonlinear systems subject to matched trigonometric-polynomial disturbances with unknown frequencies.
We first proposed a novel canonical nonlinear internal model, based on which we developed an adaptive disturbance observer that estimates such disturbances without prior knowledge of the number of sinusoidal modes or their frequencies.
The proposed design guarantees global asymptotic convergence of the disturbance estimation error without imposing any PE condition.
Moreover, under a mild PE condition, both the disturbance and parameter estimation errors converge exponentially.
Using the disturbance estimate as a feedforward input in the control design, we established sufficient conditions for tracking and disturbance rejection in general nonlinear systems, dealing with both matched and unmatched disturbances.
The effectiveness of the method was demonstrated on a trajectory tracking and disturbance rejection problem for a flexible-joint robot manipulator.
Future work will focus on extending the adaptive observer to handle the same class of disturbances in nonlinear systems with uncertain parameters by employing the internal-model reformulation technique in \cite{Wu_2025}.

% \begin{ack}
% Place acknowledgments here.
% \end{ack}

\section*{DECLARATION OF GENERATIVE AI AND AI-ASSISTED TECHNOLOGIES IN THE WRITING PROCESS}
During the preparation of this work the authors used ChatGPT in order to check grammar. After using this tool/service, the authors reviewed and edited the content as needed and take full responsibility for the content of the publication.

% \bibliography{ifacconf}             % bib file to produce the bibliography
% with bibtex (preferred)

% \bibliographystyle{plainnat}  % or another numeric natbib style                          

% \bibliographystyle{plain}        % Include this if you use bibtex
% \bibliography{reference}           % and a bib file to produce the
% bibliography (preferred). The
% correct style is generated by
% Elsevier at the time of printing.

% \begin{thebibliography}{99}     % Otherwise use the
%    % thebibliography environment.
%    % Insert the full references here.
%    % See a recent issue of Automatica
%    % for the style.

% \end{thebibliography}

% \appendix
% \section{A summary of Latin grammar}    % Each appendix must have a short title.
% \section{Some Latin vocabulary}         % Sections and subsections are supported
% % in the appendices.

\end{document}